# PREHOMOGENEOUS VECTOR SPACES AND FIELD EXTENSIONS II

Anthony C. Kable and Akihiko Yukie[1]

Oklahoma State University

**Introduction**

We fix an infinite field $k$ of any characteristic. In this paper we consider the following three prehomogeneous vector spaces

(1) $G = \mathrm{GL}(2)_{k_1} \times \mathrm{GL}(2)_k$ where $k_1/k$ is a fixed quadratic extension and $V$ is the space of pairs of binary Hermitian forms,

(2) $G = \mathrm{GL}(1)_k \times \mathrm{GL}(2)_{k_1}$ where $k_1/k$ is a fixed cubic extension and $V$ is an eight dimensional representation of $G$ which becomes the $D_4$ case in [7] after a suitable field extension of $k$,

(3) $G = \mathrm{GL}(3)_{k_1} \times \mathrm{GL}(2)_k$ where $k_1/k$ is a fixed quadratic extension and $V$ is the space of pairs of ternary Hermitian forms.

For $x \in V$ let $G_x$ be the stabilizer. For any algebraic group $G$ over $k$ we denote the connected component of 1 by $G^0$. Let $R$ be any $k$–algebra. We denote the set of invertible elements of $R$ by $R^\times$. For any variety $X$ over $k$ the set of $R$–rational points of $X$ is denoted by $X_R$. If $k^{\mathrm{sep}} \supset k' \supset k'' \supset k$ are fields such that $[k':k] < \infty$ we denote the norm and the trace by $\mathrm{N}_{k'/k''}$, $\mathrm{Tr}_{k'/k''}$.

Let $\mathfrak{S}_i$ be the permutation group of $i$ elements. As in [7] we use the notation $\mathfrak{E}\mathfrak{x}_i$ for the set of isomorphism classes of Galois extensions of $k$ which are splitting fields of degree $i$ equations without multiple roots. Note that $\mathrm{H}^1(k, \mathfrak{S}_i)$ is the set of conjugacy classes of homomorphisms from $\mathrm{Gal}(k^{\mathrm{sep}}/k)$ to $\mathfrak{S}_i$. If $i = 2$ or 3 there is a bijection between $\mathrm{H}^1(k, \mathfrak{S}_i)$ (see §1 for the definition) and $\mathfrak{E}\mathfrak{x}_i$ (see [7]). In [7] D. Wright and the second author considered eight prehomogeneous vector spaces $(G, V)$ and proved that there is a bijective correspondence between $G_k \backslash V_k^{\mathrm{ss}}$ and $\mathrm{H}^1(k, \mathfrak{S}_i)$ for suitable $i$ between 2 and 5.

The purpose of this paper is to prove an analogous correspondence for the above prehomogeneous vector spaces (1)–(3). For case (1) the correspondence is bijective. However, it turns out that the correspondence is not bijective for cases (2), (3). We describe the fiber structure of this correspondence in §§2–4. In §1 we briefly review basic properties of the non-abelian Galois cohomology set and its relation to the orbit space $G_k \setminus V_k^{\mathrm{ss}}$. In §§2–4 we consider the prehomogeneous vector spaces (1)–(3) respectively.

[1] The second author is partially supported by NSF grant DMS-9401391 and NSA grant MDA-904-93-H-3035

In §§2–4 we also determine the structure of $G_x^0$ for all $x \in V_k^{\mathrm{ss}}$ for the prehomogeneous vector spaces (1)–(3). If $k$ is a number field we can associate the zeta function for each case. The zeta function is a counting function of $G_k \setminus V_k^{\mathrm{ss}}$ possibly with the weight $\mathrm{vol}(G_{x\mathbb{A}}^0/G_{xk}^0)$. So if we determine the structure of $G_x^0$ for all $x \in V_k^{\mathrm{ss}}$, we know what kind of density theorem we can expect for each case. We discuss this issue in §5.

## §1. Rational orbits and the Galois cohomology

In this section we briefly recall the relation between the Galois cohomology set and the set of rational orbits in prehomogeneous vector spaces. Also we prove a few lemmas which we will need in later sections.

We first recall the definition of the Galois cohomology set. Let $G$ be an algebraic group over $k$, and $k'/k$ a finite Galois extension. A 1–cocycle is a function $h = \{h_\eta\}$ from $\mathrm{Gal}(k'/k)$ to $G_{k'}$ ($h_\eta$ is the value of $h$ at $\eta \in \mathrm{Gal}(k'/k)$) satisfying the condition

$$h_{\eta_1 \eta_2} = h_{\eta_2} h_{\eta_1}^{\eta_2}$$

for all $\eta_1, \eta_2$. If $h = \{h_\eta\}$, $i = \{i_\eta\}$ are 1–cocycles, they are equivalent if there exists $g \in G_{k'}$ such that

$$h_\eta = g^{-1} i_\eta g^\eta$$

for all $\eta$. This defines an equivalence relation and $\mathrm{H}^1(k'/k, G)$ is the set of equivalence classes. Let $g \in G_{k'}$. We use the notation $\delta g$ for the 1–cocycle $h = \{h_\eta\}$ defined by $h_\eta = g^{-1} g^\eta$ for all $\eta \in \mathrm{Gal}(k'/k)$. The cohomology class defined by $\delta g$ does not depend on the choice of $g$ and we denote this element by 1.

We define $\mathrm{H}^1(k, G)$ to be the projective limit of $\mathrm{H}^1(k'/k, G)$ for all the finite Galois extensions $k'$. We define $\mathrm{H}^0(k'/k, G) = \mathrm{H}^0(k, G) = G_k$. If $G$ is an abelian group, $\mathrm{H}^n(k'/k, G)$ can be defined for all $n$ and has a structure of an abelian group also.

Let

(1.1) $$1 \to G_1 \to G_2 \to G_3 \to 1$$

be a short exact sequence of algebraic groups over $k$. This means that $G_1$ is a normal subgroup of $G_2$, the kernel of $G_2 \to G_3$ is $G_1$, and $G_{2\,k^{\mathrm{sep}}} \to G_{3\,k^{\mathrm{sep}}}$ is surjective. If $G_1, G_2, G_3$ are abelian groups, we have the following long exact sequence

$$\cdots \to \mathrm{H}^n(k, G_1) \to \mathrm{H}^n(k, G_2) \to \mathrm{H}^n(k, G_3) \to \cdots.$$

We consider the case when $G_1, G_2, G_3$ are not necessarily abelian. Let $g \in G_{3\,k}$. If $k'/k$ is a large enough finite Galois extension, there is an element $f \in G_{2\,k'}$ which maps to $g$. For a cohomology class $c$ in $\mathrm{H}^1(k, G_1)$ defined by a 1–cocycle $h = \{h_\eta\}$, we define $gc \in \mathrm{H}^1(k, G_1)$ to be the class defined by the 1–cocycle $\{f h_\eta (f^\eta)^{-1}\}$. Since $g \in G_{3\,k}$, $f h_\eta (f^\eta)^{-1} \in G_{1\,k^{\mathrm{sep}}}$ for all $\eta \in \mathrm{Gal}(k^{\mathrm{sep}}/k)$ and it is easy to see that $gc$ does not depend on the choice of $f$ or $k'$. This defines an action of $G_{3\,k}$ on $\mathrm{H}^1(k, G_1)$. The following lemma is an easy consequence of Proposition 38, §5.5 and Corollaire 1, §5.5 of [6].

**Lemma (1.2)** *The sequence*

$$1 \to G_{3\,k} \backslash \mathrm{H}^1(k, G_1) \to \mathrm{H}^1(k, G_2) \to \mathrm{H}^1(k, G_3)$$

*is exact. Moreover, if* (1.1) *is split, the last map is surjective.*

Note that the exactness of the sequence in (1.2) means that the inverse image of $1 \in \mathrm{H}^1(k, G_3)$ is $G_{3\,k} \backslash \mathrm{H}^1(k, G_1)$.

For the prehomogeneous vector spaces (1)–(3), we prove in §§2–4 that there is a distinguished element $w \in V_k^{\mathrm{ss}}$ and a split exact sequence

(1.3) $$1 \to G_w^0 \to G_w \to \mathfrak{S}_i \to 1$$

where the Galois group acts on $\mathfrak{S}_i$ trivially.

It is a familiar fact that both $\mathrm{H}^1(k, \mathrm{GL}(n))$ and $\mathrm{H}^1(k, \mathrm{SL}(n))$ are trivial. As remarked in [6] (see the proof of Théorème 1, §2.2) one has

$$\mathrm{H}^1(k_1, G) = \mathrm{H}^1(k, R_{k_1/k}(G))$$

for any algebraic group $G$ over $k_1$, where $R_{k_1/k}$ denotes restriction of scalars. This leads at once to the following.

**Lemma (1.4)** *Let $k_1/k$ be a finite separable extension, and $G = \mathrm{GL}(n)_{k_1}$ or $\mathrm{SL}(n)_{k_1}$ considered as an algebraic group over $k$. Then $\mathrm{H}^1(k, G) = \{1\}$.*

In §§2–4 we prove that the prehomogeneous vector spaces (1)–(3) satisfy the following condition.

**Condition (1.5)** *The set $V_{k^{\mathrm{sep}}}^{\mathrm{ss}}$ is a single $G_{k^{\mathrm{sep}}}$–orbit.*

Suppose Condition (1.5) is satisfied. Then for any $x \in V_k^{\mathrm{ss}}$, we can choose a finite Galois extension $k'/k$ and $g \in G_{k'}$ such that $x = gw$. Then $c_x = \{g^{-1}g^\eta\}$ determines an element of $\mathrm{Ker}(\mathrm{H}^1(k, G_w) \to \mathrm{H}^1(k, G))$ (which is the set of elements which map to $1 \in \mathrm{H}^1(k, G)$). In [2] Igusa assumed that the characteristic of the field is zero. However, if Condition (1.5) is satisfied, we can still make cohomology classes from rational orbits in $V_k^{\mathrm{ss}}$. Therefore, without changing Igusa's argument, we have the following Theorem.

**Theorem (1.6) (Igusa)** *Suppose a prehomogeneous vector space $(G, V)$ satisfies Condition (1.5). Then the correspondence*

$$G_k \backslash V_x^{\mathrm{ss}} \ni x \to c_x \in \mathrm{Ker}(\mathrm{H}^1(k, G_w) \to \mathrm{H}^1(k, G))$$

*is bijective.*

For the prehomogeneous vector spaces (1)–(3), $G$ is either $\mathrm{GL}(2)_{k_1} \times \mathrm{GL}(2)_k$, $\mathrm{GL}(1)_k \times \mathrm{GL}(2)_{k_1}$, or $\mathrm{GL}(2)_{k_1} \times \mathrm{GL}(3)_k$, where $k_1$ is either a quadratic or cubic extension of $k$. So by Lemma (1.4) and Theorem (1.6), we have the following proposition.

**Proposition (1.7)** *For the prehomogeneous vector spaces (1)–(3), the correspondence*
$$G_k \backslash V_k^{\mathrm{ss}} \ni x \to c_x \in \mathrm{H}^1(k, G_w)$$
*is bijective.*

Since (1.3) is a split exact sequence, by Lemma (1.2), we have the following exact sequence

(1.8) $$1 \to \mathfrak{S}_i \backslash \mathrm{H}^1(k, G_w^0) \to \mathrm{H}^1(k, G_w) \to \mathrm{H}^1(k, \mathfrak{S}_i) \to 1.$$

Therefore, the canonical map $\mathrm{H}^1(k, G_w) \to \mathrm{H}^1(k, \mathfrak{S}_i)$ can be considered as a map from $\mathrm{H}^1(k, G_w) \cong G_k \backslash V_k^{\mathrm{ss}}$ to $\mathfrak{E}\mathfrak{x}_i$. We denote this map by $\alpha_V$. If $x \in V_k^{\mathrm{ss}}$, we also use the notation $\alpha_V(x)$ for $\alpha_V(G_k x)$ and call this field $k(x)$. In §§2,4 we define a subscheme $\mathrm{Zero}(x) \subset \mathbb{P}^1$ defined over $k$ for any $x \in V_k^{\mathrm{ss}}$. It has the property that $k(x)$ coincides with the field generated by residue fields of points in $\mathrm{Zero}(x)$. Moreover, from the naturality of the construction of $\mathrm{Zero}(x)$, it will turn out that the following sequence

(1.9) $$1 \to G_x^0 \to G_x \to \mathrm{Aut}(\mathrm{Zero}(x)) \to 1$$

is exact (but not necessarily split). Here $\mathrm{Aut}(\mathrm{Zero}(x))$ is the algebraic group over $k$ which represents the functor $S \to \mathrm{Aut}_S(\mathrm{Zero}(x) \times_k S)$ for $k$–schemes $S$. In §3 we prove that for any $k' \in \mathfrak{E}\mathfrak{x}_2$, there is an element $x \in \alpha_V^{-1}(k')$ such that there is a split exact sequence

(1.10) $$1 \to G_x^0 \to G_x \to \mathfrak{S}_2 \to 1,$$

where the Galois group acts trivially on $\mathfrak{S}_2$.

Let $x \in V_k^{\mathrm{ss}}$. We choose an element $g_x \in G_{k^{\mathrm{sep}}}$ so that $x = g_x w$. Then for each element $c \in \mathrm{H}^1(k, G_w)$ defined by a 1–cocycle $\{h_\eta\}$, we can associate an element $c^{g_x} \in \mathrm{H}^1(k, G_x)$ defined by a 1–cocycle $\{g_x h_\eta (g_x^\eta)^{-1}\}$. It is easy to see that the map $c \to c^{g_x}$ is well defined and does not depend on the choice of $g_x$. Also a similar construction using $g_x^{-1}$ defines a map from $\mathrm{H}^1(k, G_x)$ to $\mathrm{H}^1(k, G_w)$. Therefore, we have the following lemma.

**Lemma (1.11)** *The map*
$$\mathrm{H}^1(k, G_w) \ni c \to c^{g_x} \in \mathrm{H}^1(k, G_x)$$
*induces a bijection.*

In the following lemma let $i = 2$ for cases (1), (2) and $i = 3$ for case (3).
Consider $x \in V_k^{\mathrm{ss}}$ in (1.9) or (1.10).

**Lemma (1.12)** *Let $k' = k(x) \in \mathfrak{E}\mathfrak{x}_i$. Then*
$$\alpha_V^{-1}(k') \cong \mathrm{Aut}_k(\mathrm{Zero}(x)) \backslash \mathrm{H}^1(k, G_x^0) \text{ or } \mathfrak{S}_2 \backslash \mathrm{H}^1(k, G_x^0).$$

*Moreover, by this identification, the cohomology class $\{g^{-1} g^\eta\} \in \mathrm{H}^1(k, G_x^0)$ ($g \in G_{k^{\mathrm{sep}}}$) corresponds to the orbit $G_k g x$.*

*Proof.* Let $x = g_x w$ and $y \in V_k^{\text{ss}}$. Then $y = gx$ for certain $g \in G_{k^{\text{sep}}}$. So $y = gg_x w$, which implies that $y$ corresponds to the cohomology class $c = \{(gg_x)^{-1}(gg_x)^\eta\}$. Therefore, by the identification in Lemma (1.11),

$$c^{g_x} = \{g_x(gg_x)^{-1}(gg_x)^\eta(g_x^\eta)^{-1}\} = \{g^{-1}g^\eta\}.$$

So if we identify $G_x \setminus V_k^{\text{ss}}$ with $\text{H}^1(k, G_x)$, the orbit of $y$ corresponds to the cohomology class $c_{x,y} = \{g^{-1}g^\eta\}$. Because of Lemma (1.2), we only have to prove that $\alpha_V(y) = \alpha_V(x)$ if and only if $c_{x,y}$ comes from $\text{H}^1(k, G_x^0)$.

Since both $x$ and $w$ are rational elements, $x = g_x^\eta w$ also. Therefore,

$$G_{xk^{\text{sep}}} = g_x G_{wk^{\text{sep}}} g_x^{-1} = g_x^\eta G_{wk^{\text{sep}}}(g_x^\eta)^{-1},$$
$$G^0_{xk^{\text{sep}}} = g_x G^0_{wk^{\text{sep}}} g_x^{-1} = g_x^\eta G^0_{wk^{\text{sep}}}(g_x^\eta)^{-1}.$$

Since (1.3) is split, we consider $\mathfrak{S}_i$ as a subgroup of $G_w$. Suppose $y = gx = gg_x w$ satisfies the condition $\alpha_V(y) = \alpha_V(x)$. Then $\{g_x^{-1}g^{-1}g^\eta g_x^\eta\}$ and $\{g_x^{-1}g_x^\eta\}$ map to the same element in $\text{H}^1(k, \mathfrak{S}_i)$. So there exist $r \in \mathfrak{S}_i$ and $n_\eta \in G^0_{wk^{\text{sep}}}$ such that

$$g_x^{-1}g^{-1}g^\eta g_x^\eta = r^{-1}g_x^{-1}g_x^\eta r n_\eta.$$

We can modify the above equation as

$$(g_x r^{-1} g_x^{-1})^{-1}(g^{-1}g^\eta)(g_x r^{-1} g_x^{-1})^\eta = g_x^\eta r n_\eta r^{-1}(g_x^\eta)^{-1}.$$

The left hand side defines the same cohomology class as $\{g^{-1}g^\eta\}$ in $\text{H}^1(k, G_x)$ and the right hand side belongs to $G^0_{xk^{\text{sep}}}$. Therefore, this cohomology class comes from $\text{H}^1(k, G_x^0)$.

Conversely, if $\{g^{-1}g^\eta\}$ comes from $\text{H}^1(k, G_x^0)$, by changing $g$ is necessary, we can assume that $g^{-1}g^\eta \in G^0_{xk^{\text{sep}}}$ for all $\eta$. Then

$$g_x^{-1}g^{-1}g^\eta g_x^\eta = g_x^{-1}g_x^\eta(g_x^\eta)^{-1}g^{-1}g^\eta g_x^\eta,$$

and

$$(g_x^\eta)^{-1}g^{-1}g^\eta g_x^\eta \in G^0_{wk^{\text{sep}}}$$

for all $\eta$. Therefore, $\alpha_V(y) = \alpha_V(x)$. This proves the lemma.

Q.E.D.

**Remark (1.13)** In [7] the logic was slightly imprecise. In order to determine the fiber structure of $\alpha_V$, Lemma (1.12) should have been used. However, since $\text{H}^1(k, G_x^0) = 1$ for all the cases in [7] by Lemma (1.4), the statements in [7] do not have to be changed. Also, in order to apply Igusa's result, we have to show that Condition (1.5) is satisfied. For the $F_4$ and $E_8$ cases in [7], it follows from Propositions 1.1, 1.5, 1.6, 2.3, Lemma 2.10, and Corollary 2.12 (all in [7]). Other cases are straightforward. Also the assumption in [7] that the characteristic of the field is not 2,3, or 5 was not necessary. This is because if the discriminant of a rational polynomial is not zero, its roots generate a separable extension no matter what the characteristic of the field is. For example, the discriminant of the quadratic polynomial $av_1^2 + bv_1v_2 + cv_2^2$ is $b^2 - 4ac$. If the characteristic of the field is 2, it is $b^2$. If $b \neq 0$, this polynomial cannot be of the form $(\alpha v_1 + \beta v_2)^2 = \alpha^2 v_1^2 + \beta^2 v_2^2$.

### §2. The non-split $D_4$ case (1)

In this section, we consider the space of pairs of binary Hermitian forms and prove that the set $G_k \backslash V_k^{\text{ss}}$ corresponds bijectively with $\mathfrak{Ex}_2$. We also determine the stabilizer of any element $x \in V_k^{\text{ss}}$.

We fix a separable quadratic extension $k_1 = k(\alpha_0)$ of $k$. The non-trivial element of $\text{Gal}(k_1/k)$ is denoted by $\sigma$. Let $G = \text{GL}(2)_{k_1} \times \text{GL}(2)_k$ considered as an algebraic group over $k$. Let $W$ be the space of binary Hermitian forms, i.e. any element in $W$ is a $2 \times 2$ matrix $A$ satisfying ${}^t A^\sigma = A$. Let $V = W \otimes k^2$. We consider $V$ as the space of binary Hermitian forms $M(v)$ with entries in the space of linear forms in two variables $v = (v_1, v_2)$. Then $g = (g_1, g_2) \in G$ acts on $V$ by $gM(v) = g_1 M(vg_2){}^t g_1^\sigma$.

There is a natural map $W \otimes k_1^2 \to \text{M}(2,2)_{k_1} \otimes k_1^2$. This map is equivariant with respect to $G_{k_1} \cong \text{GL}(2)_{k_1} \times \text{GL}(2)_{k_1} \times \text{GL}(2)_{k_1}$. Since $\text{M}(2,2)_{k_1} \otimes k_1^2$ is an irreducible representation of $G_{k_1}$, this map is surjective. Since the dimension of $W \otimes k_1^2$ is eight, this map is an isomorphism. Therefore, $(G, V)$ is a prehomogeneous vector space and is a $k$–form of the $D_4$ case in [7]. Note that this argument works even if the characteristic of the field is two.

Let $x = v_1 x_1 + v_2 x_2$ where $x_1, x_2$ are binary Hermitian matrices. Consider the map
$$x \to F_x(v) = \det x \in \text{Sym}^2 k^2.$$
Note that since $x_1, x_2$ are Hermitian, $\det x \in \text{Sym}^2 k^2$. Clearly, $F_{(g_1,g_2)x}(v) = \text{N}_{k_1/k}(\det g_1) F_x(vg_2)$. Let $\text{Zero}(x) \subset \mathbb{P}^1$ be the subset defined by the roots of $F_x(v) = 0$. More precisely,
$$\text{Zero}(x) = \text{Proj } k[v_1, v_2]/(F_x(v))$$
as a scheme over $k$. Let $\Delta(x)$ be the discriminant of $F_x(v)$ as a polynomial of $v$. Clearly $\Delta$ is a non-constant relative invariant polynomial. Since $\Delta$ does not vanish, $\text{Zero}(x)$ is a reduced scheme for every $x \in V_k^{\text{ss}}$.

It is possible to check by linear algebra that the $D_4$ case in [7] satisfies Condition (1.5). Since $k_1/k$ is a separable extension, we get the following proposition.

**Proposition (2.1)** *The prehomogeneous vector space $(G, V)$ satisfies* Condition (1.5).

If $g = (g_1, g_2) \in G_x$, $v \to v g_1^{-1}$ is an automorphism of $\text{Zero}(x)$. So there is a natural homomorphism $G_x \to \text{Aut}(\text{Zero}(x))$.

Let

(2.2)
$$w = v_1 \begin{pmatrix} 1 & \\ & 0 \end{pmatrix} + v_2 \begin{pmatrix} 0 & \\ & 1 \end{pmatrix},$$
$$\tau = \left( \begin{pmatrix} 0 & 1 \\ 1 & 0 \end{pmatrix}, \begin{pmatrix} 0 & 1 \\ 1 & 0 \end{pmatrix} \right).$$

Note that $\tau \in G_{wk}$. As in [7], $\text{Zero}(w) = \{(1,0), (0,1)\}$, $\text{Aut}(\text{Zero}(x)) \cong \mathfrak{S}_2$, and $\tau$ exchanges $(1,0)$ and $(0,1)$.

For the rest of this paper we use the notation

(2.3) $$a(t_1, t_2) = \begin{pmatrix} t_1 & 0 \\ 0 & t_2 \end{pmatrix}, \quad n(u) = \begin{pmatrix} 1 & 0 \\ u & 1 \end{pmatrix}.$$

Let

(2.4) $$t = (a_2(t_{11}, t_{12}), a_2(t_{21}, t_{22})),$$

where $t_{11}, t_{12} \in k_1^\times$, $t_{21}, t_{22} \in k^\times$. By considering the stabilizer of $w$ over $k^{\text{sep}}$ we find that $\text{Ker}(G_w \to \mathfrak{S}_2) = G_w^0$ and that if $g = (g_1, g_2) \in G_w^0$ then $g$ must have the form (2.4). Since $\tau \in G_{wk}$ we have a split exact sequence

(2.5) $$1 \to G_w^0 \to G_w \to \mathfrak{S}_2 \to 1,$$

where the action of the Galois group on $\mathfrak{S}_2$ is trivial. Now $t$ in (2.4) belongs to $G_w$ if and only if

$$t_{21} \text{N}_{k_1/k}(t_{11}) = t_{22} \text{N}_{k_1/k}(t_{12}) = 1.$$

Therefore, we get the following proposition.

**Proposition (2.6)** *As an algebraic group over $k$,*

$$G_w^0 = \{t \mid t_{21} \text{N}_{k_1/k}(t_{11}) = t_{22} \text{N}_{k_1/k}(t_{12}) = 1\}$$
$$\cong \text{GL}(1)_{k_1} \times \text{GL}(1)_{k_1}.$$

In general, if $x \in V_k^{\text{ss}}$ and $x = gw$ for $g = (g_1, g_2) \in G_{k^{\text{sep}}}$,

(2.7) $$G_{xk^{\text{sep}}}^0 = g G_{wk^{\text{sep}}}^0 g^{-1}, \quad G_{xk^{\text{sep}}} = g G_{wk^{\text{sep}}} g^{-1},$$
$$\text{Zero}(x)_{k^{\text{sep}}} = \{(1, 0) g_1^{-1}, (0, 1) g_1^{-1}\}.$$

Therefore the sequence (1.9) is exact.

By (1.8) we have a surjective map $\alpha_V : G_k \backslash V_k^{\text{ss}} \to \mathfrak{E}\mathfrak{x}_2$. Let $k(\alpha)/k$ be a separable quadratic extension. If $\alpha \notin k_1$,

$$G_{k(\alpha)} \cong \text{GL}(2)_{k_1(\alpha)} \times \text{GL}(2)_{k(\alpha)},$$
$$V_{k(\alpha)} \cong W_{k(\alpha)} \otimes k(\alpha)^2.$$

If $\alpha = \alpha_0$,

$$G_{k_1} \cong \text{GL}(2)_{k_1} \times \text{GL}(2)_{k_1} \times \text{GL}(2)_{k_1},$$
$$V_{k_1} \cong \text{M}(2,2)_{k_1} \otimes k_1^2.$$

Let $g = (g_1, g_2) \in G_k$ and $x = v_1 x_1 + v_2 x_2$. Under the above identification $g$ corresponds to $(g_1, g_1^\sigma, g_2)$. The action of the Galois group $\text{Gal}(k_1/k)$ on $G_{k_1}$ is given by $g^\sigma = (g_2^\sigma, g_1^\sigma, g_3^\sigma)$ ($\sigma \in \text{Gal}(k_1/k)$ is the non-trivial element). There is a unique right action of the Galois group $\text{Gal}(k_1/k)$ on $W \otimes k_1$ satisfying the property that it is trivial on $W$ and $(tx)\eta = t^\eta x\eta$ for $\eta \in \text{Gal}(k_1/k)$. Therefore, $x \mapsto {}^t x^\sigma$ for the non-trivial element $\sigma \in \text{Gal}(k_1/k)$, is the Galois group action on $\text{M}(2,2)_{k_1}$ induced by that on $W$. This implies that the action of the Galois group $\text{Gal}(k_1/k)$ on $V_{k_1}$ is given by $x^\sigma = v_1{}^t x_1^\sigma + v_2{}^t x_2^\sigma$. Also the action of $G_{k_1}$ on $V_{k_1}$ is given by $gM(v) = g_1 M(vg_3){}^t g_2$ for $g = (g_1, g_2, g_3) \in G_{k_1}$.

Let $f(v) = v_1^2 + a_1 v_1 v_2 + a_2 v_2^2 \in k[v]$ be a polynomial such that the roots $\alpha_1, \alpha_2$ of $f$ generate the field $k(\alpha)$. (This implies that $\alpha_1 + \alpha_2 = -a_1$, $\alpha_1 \alpha_2 = a_2$.)

We define

$$(2.8) \quad g_\alpha = \begin{cases} \left( \begin{pmatrix} 1 & 1 \\ -\alpha_1 & -\alpha_2 \end{pmatrix}, \begin{pmatrix} 1 & 1 \\ -\alpha_1 & -\alpha_2 \end{pmatrix} \right) & \alpha \notin k_1, \\ \left( \begin{pmatrix} 1 & 1 \\ -\alpha_1 & -\alpha_2 \end{pmatrix}, \begin{pmatrix} 1 & 1 \\ -\alpha_1 & -\alpha_2 \end{pmatrix}, \begin{pmatrix} 1 & 1 \\ -\alpha_1 & -\alpha_2 \end{pmatrix} \right) & \alpha = \alpha_0, \end{cases}$$

$$w_\alpha = v_1 \begin{pmatrix} 2 & a_1 \\ a_1 & a_1^2 - 2a_2 \end{pmatrix} + v_2 \begin{pmatrix} a_1 & a_1^2 - 2a_2 \\ a_1^2 - 2a_2 & a_1^3 - 3a_1 a_2 \end{pmatrix}.$$

Then $w_\alpha = g_\alpha w \in V_k^{\text{ss}}$.

It is easy to see that

$$F_{w_\alpha}(v) = (\alpha_1 - \alpha_2)^2 (v_1 - \alpha_1 v_2)(v_1 - \alpha_2 v_2) = (\alpha_1 - \alpha_2)^2 f(v).$$

The field $k(w_\alpha)$ corresponds by definition to the cohomology class $\{g_\alpha^{-1} g_\alpha^\eta\} \in \mathrm{H}^1(k, \mathfrak{S}_2)$. If $\nu$ is the non-trivial element of $\mathrm{Gal}(k(\alpha)/k)$ then $g_\alpha^{-1} g_\alpha^\nu = \tau$ and so this cohomology class also corresponds to the field $k(\alpha)$. Therefore $k(w_\alpha) = k(\alpha)$.

Let

$$(2.9) \quad t = \begin{cases} (a_2(t_{11}, t_{12}), a_2(t_{21}, t_{22})) & \alpha \notin k_1, \\ (a_2(t_{11}, t_{12}), a_2(t_{21}, t_{22}), a_2(t_{31}, t_{32})) & \alpha = \alpha_0. \end{cases}$$

**Proposition (2.10)** (1) *If $\alpha \notin k_1$, as an algebraic group over $k$,*

$$G_{w_\alpha}^0 = \left\{ g_\alpha t g_\alpha^{-1} \;\middle|\; \begin{array}{l} t_{11}, t_{12} \in k_1(\alpha)^\times,\ t_{21}, t_{22} \in k(\alpha)^\times,\ t_{12} = t_{11}^\nu, \\ t_{21} \mathrm{N}_{k_1(\alpha)/k(\alpha)}(t_{11}) = t_{22} \mathrm{N}_{k_1(\alpha)/k(\alpha)}(t_{12}) = 1 \end{array} \right\}$$
$$\cong \mathrm{GL}(1)_{k_1(\alpha)}.$$

(2) *If $\alpha = \alpha_0$, as an algebraic group over $k$,*

$$G_{w_\alpha}^0 = \{ g_\alpha t g_\alpha^{-1} \mid t_{ij} \in k_1^\times \text{ for all } i,j,\ t_{32} = t_{31}^\sigma,\ t_{21} = t_{12}^\sigma,\ t_{22} = t_{11}^\sigma,\ t_{31} t_{11} t_{12}^\sigma = 1 \}$$
$$\cong \mathrm{GL}(1)_{k_1} \times \mathrm{GL}(1)_{k_1}.$$

*Proof.* Consider (1). In order to prove an isomorphism between two algebraic groups $G_1, G_2$ over $k$, it is enough to prove natural isomorphisms between the sets $G_{1R}, G_{2R}$ of $R$–rational points of $G_1, G_2$ for all $k$–algebras $R$. For this, the reader should see Theorem [3, p. 17].

Let $R$ be any $k$–algebra. For any Galois extension $k'/k$, $\nu \in \mathrm{Gal}(k'/k)$ acts on $k' \otimes R$ by $(x \otimes r)^\nu = x^\nu \otimes r$. We define $R(\alpha) = k(\alpha) \otimes R$ and $R_1(\alpha) = k_1(\alpha) \otimes R$. Then

$$G_{w_\alpha R}^0 = \{ g \in G_{w_\alpha R(\alpha)}^0 \mid g^\nu = g \text{ for all } \nu \in \mathrm{Gal}(k(\alpha)/k) \}.$$

Any element of $G^0_{w_\alpha R(\alpha)}$ is of the form $g_\alpha t g_\alpha^{-1}$ where $t$ is as in (2.4) with $t_{11}, t_{12} \in R_1(\alpha)^\times$ and $t_{21}, t_{22} \in R(\alpha)^\times$. Since $g_\alpha^\nu = g_\alpha \tau$,

$$(g_\alpha t g_\alpha^{-1})^\nu = g_\alpha \tau (a_2(t_{11}^\nu, t_{12}^\nu), a_2(t_{21}^\nu, t_{22}^\nu)) \tau g_\alpha^{-1}$$
$$= g_\alpha (a_2(t_{12}^\nu, t_{11}^\nu), a_2(t_{22}^\nu, t_{21}^\nu)) g_\alpha^{-1}.$$

So $(g_\alpha t g_\alpha^{-1})^\nu \in G^0_{w_\alpha R}$ if and only if

$$t_{12} = t_{11}^\nu, \ t_{21} N_{k_1/k}(t_{11}) = t_{22} N_{k_1/k}(t_{12}) = 1.$$

Note that this condition implies $t_{22} = t_{21}^\nu$ also. This proves the statement (1).

The statement (2) follows by a similar argument using the fact

$$(g_\alpha t g_\alpha^{-1})^\sigma = g_\alpha \tau (a_2(t_{21}^\sigma, t_{22}^\sigma), a_2(t_{11}^\sigma, t_{12}^\sigma), a_2(t_{31}^\sigma, t_{32}^\sigma)) \tau g_\alpha^{-1}$$
$$= g_\alpha (a_2(t_{22}^\sigma, t_{21}^\sigma), a_2(t_{12}^\sigma, t_{11}^\sigma), a_2(t_{32}^\sigma, t_{31}^\sigma)) g_\alpha^{-1}.$$

Q.E.D.

By Lemma (1.4) and Proposition (2.10), $H^1(k, G^0_{w_\alpha}) = \{1\}$ for all the cases. Therefore, we have the following theorem.

**Theorem (2.11)** *The map $\alpha_V : G_k \backslash V_k^{ss} \to \mathfrak{E}\mathfrak{x}_2$ is bijective.*

For $x = w_\alpha$, the field $k(x)$ is generated by residue fields of points in $\mathrm{Zero}(x)$. But by the above theorem, all the points in $V_k^{ss}$ are either in $G_k w$ or in $G_k w_\alpha$ for some $\alpha$. Therefore, we get the following corollary.

**Corollary (2.12)** *If $x \in V_k^{ss}$, the field $k(x)$ is generated by residue fields of points in $\mathrm{Zero}(x)$.*

### §3. The non-split $D_4$ case (2)

In this section, we consider the most non-split prehomogeneous vector space which becomes the $D_4$ case in [7] after a suitable extension of the base field. We first describe the prehomogeneous vector space we consider in this section.

Let $k_1$ be a separable cubic extension of $k$, and $k_2$ the normal closure of $k_1$. Then either $k_2 = k_1$ is a cyclic cubic extension of $k$ or $k_2$ is an $\mathfrak{S}_3$–extension of $k$. Let $G = \mathrm{GL}(2)_{k_1}$ considered as a group over $k$. We construct a prehomogeneous vector space $(G, V)$ over $k$ such that $(G_{k_2}, V_{k_2})$ is the prehomogeneous vector space $(\widetilde{G}, \widetilde{V})$ where

$$\widetilde{G} = \mathrm{GL}(2)_{k_2} \times \mathrm{GL}(2)_{k_2} \times \mathrm{GL}(2)_{k_2}, \ \widetilde{V} = k_2^2 \otimes k_2^2 \otimes k_2^2.$$

We choose three different imbeddings $\sigma_1, \sigma_2, \sigma_3 : k_1 \to k_2$ over $k$. We consider $G_{k_2}$ as a subset of $\mathrm{M}(2,2)_{k_1} \otimes k_2$. Then the map

$$\phi : \mathrm{M}(2,2)_{k_1} \otimes k_2 \ni g \otimes a \to (ag^{\sigma_1}, ag^{\sigma_2}, ag^{\sigma_3}) \in \mathrm{M}(2,2)_{k_2} \times \mathrm{M}(2,2)_{k_2} \times \mathrm{M}(2,2)_{k_2}$$

induces an isomorphism. By this map, $\phi(G_{k_2}) = \mathrm{GL}(2)_{k_2} \times \mathrm{GL}(2)_{k_2} \times \mathrm{GL}(2)_{k_2}$.

Let $H_1 = \mathrm{Gal}(k_2/k)$ and $H_2 = \mathrm{Gal}(k_2/k_1)$. Then $H_2$ is a subgroup of $H_1$ and $[H_1 : H_2] = 3$. So there exists a homomorphism $h : H_1 \to \mathfrak{S}_3$ such that $H_2 \sigma_i \sigma = H_2 \sigma_{h(\sigma)(i)}$ for $i = 1, 2, 3$. (Here if $\tau_1, \tau_2 \in \mathfrak{S}_3$, $(\tau_1 \tau_2)(i) = \tau_2(\tau_1(i))$.) If $g \in \mathrm{GL}(2)_{k_1}$, we regard $g$ as a $k$–rational point of $G$. So $g$ is fixed by the action of $H_1$. This means that if $a \in k_2$ and $\sigma \in H_1$,

$$\begin{aligned}
\phi((g \otimes a)^\sigma) &= \phi(g \otimes a^\sigma) \\
&= (a^\sigma g^{\sigma_1}, a^\sigma g^{\sigma_2}, a^\sigma g^{\sigma_2}) \\
&= ((ag^{\sigma_1 \sigma^{-1}})^\sigma, (ag^{\sigma_2 \sigma^{-1}})^\sigma, (ag^{\sigma_2 \sigma^{-1}})^\sigma) \\
&= ((ag^{\sigma_{h(\sigma^{-1})(1)}})^\sigma, (ag^{\sigma_{h(\sigma^{-1})(2)}})^\sigma, (ag^{\sigma_{h(\sigma^{-1})(3)}})^\sigma).
\end{aligned}$$

Therefore, if $g_1, g_2, g_3 \in \mathrm{GL}(2)_{k_2}$,

(3.1) $$(g_1, g_2, g_3)^\sigma = (g^\sigma_{h(\sigma^{-1})(1)}, g^\sigma_{h(\sigma^{-1})(2)}, g^\sigma_{h(\sigma^{-1})(3)}).$$

Let $\widetilde{V} = k_2^2 \otimes k_2^2 \otimes k_2^2$. We define a right action of $H_1$ on $\widetilde{V}$ by

(3.2) $$(x_1 \otimes x_2 \otimes x_3)^\sigma = x^\sigma_{h(\sigma^{-1})(1)} \otimes x^\sigma_{h(\sigma^{-1})(2)} \otimes x^\sigma_{h(\sigma^{-1})(3)}.$$

Let $V = \widetilde{V}^{H_1}$. Clearly, $(gx)^\sigma = g^\sigma x^\sigma$ for $\sigma \in H_1$, $g \in \widetilde{G}$, and $x \in \widetilde{V}$. Since $G_k = G^{H_1}_{k_2}$, $G$ naturally acts on $V$. We show that $V \otimes k_2 \cong \widetilde{V}$.

We first fix a coordinate system for $\widetilde{V}$. Let $f_1 = \begin{pmatrix} 1 \\ 0 \end{pmatrix}$, $f_2 = \begin{pmatrix} 0 \\ 1 \end{pmatrix}$, and $e_{ijk} = f_i \otimes f_j \otimes f_k$ for $i, j, k = 1, 2$. Then $\{e_{ijk} \mid i, j, k = 1, 2\}$ is a basis for $\widetilde{V}$. So any element $x \in \widetilde{V}$ can be expressed as $x = \sum_{i,j,k} x_{ijk} e_{ijk}$, where $x_{ijk} \in k_2$ for all $i, j, k$.

For the rest of this section, we choose and fix $\sigma \in H_1$ so that $h(\sigma) = (123)$. If $H_1 \cong \mathfrak{S}_3$, without loss of generality, we may assume that $h(\bar{\sigma}) = (23)$ for the non-trivial element $\bar{\sigma} \in H_2$. Note that if $x_1, x_2, x_3 \in k$, $(x_1 \otimes x_2 \otimes x_3)^\sigma = x_3 \otimes x_1 \otimes x_2$. Therefore
$$e^\sigma_{111} = e_{111}, \ e^\sigma_{112} = e_{211}, \ e^\sigma_{121} = e_{112}, \ e^\sigma_{122} = e_{212},$$
$$e^\sigma_{211} = e_{121}, \ e^\sigma_{212} = e_{221}, \ e^\sigma_{221} = e_{122}, \ e^\sigma_{222} = e_{222}.$$

So if $x = \sum_{i,j,k} x_{ijk} e_{ijk}$,

$$\begin{aligned}
x^\sigma = &\, x^\sigma_{111} e_{111} + x^\sigma_{112} e_{211} + x^\sigma_{121} e_{112} + x^\sigma_{122} e_{212} \\
&+ x^\sigma_{211} e_{121} + x^\sigma_{212} e_{221} + x^\sigma_{221} e_{122} + x^\sigma_{222} e_{222}.
\end{aligned}$$

It is easy to see that $x^{\bar\sigma} = \sum_{i,j,k} x^{\bar\sigma}_{ikj} e_{ijk}$.

Therefore the condition $x \in V$ is equivalent to

(3.3) $x_{111}, x_{222} \in k$, $x^\sigma_{122} = x_{212}$, $x^{\sigma^2}_{122} = x_{221}$, $x^\sigma_{211} = x_{121}$, $x^{\sigma^2}_{211} = x_{112}$, $x_{122}, x_{211} \in k_1$.

So $x$ is determined by $x_{111}, x_{222} \in k$ and $x_{122}, x_{211} \in k_1$.

Since $\widetilde{V}$ is an irreducible representation of $\mathrm{GL}(2)_{k_2} \times \mathrm{GL}(2)_{k_2} \times \mathrm{GL}(2)_{k_2}$, the natural map $V \otimes k_2 \to \widetilde{V}$ is surjective. Since $\dim_k V = \dim_{k_2} \widetilde{V} = 8$, $V \otimes k_2 \cong \widetilde{V}$.

Therefore, $(G, V)$ is a prehomogeneous vector space. Since $k_2/k$ is a separable extension, we get the following proposition.

**Proposition (3.4)** *The prehomogeneous vector space $(G, V)$ satisfies* Condition (1.5).

In order to describe the action of $G$ on $V$, it is enough to consider the action of elements $a(t_1, t_2), n(u)$ and

$$(3.5) \qquad \tau = \begin{pmatrix} 0 & 1 \\ 1 & 0 \end{pmatrix}.$$

Note that this element $\tau$ is different from that in (2.2).

We choose $\sigma_1 = 1$, $\sigma_2 = \sigma$, $\sigma_3 = \sigma^2$. If $k_2$ is an $\mathfrak{S}_3$-extension of $k$, this means that we are restricting elements of $\mathrm{Gal}(k_2/k)$ to $k_1$. It is easy to see that

$$(3.6) \qquad a(t_1, t_2)x = \mathrm{N}_{k_1/k}(t_1) x_{111} e_{111} + t_1 t_2^{-1} \mathrm{N}_{k_1/k}(t_2) x_{122}$$
$$+ t_2 t_1^{-1} \mathrm{N}_{k_1/k}(t_1) x_{211} + \cdots + \mathrm{N}_{k_1/k}(t_2) x_{222} e_{222}.$$

By easy computations,

$$n(u) e_{111} = e_{111} + u e_{211} + u^\sigma e_{121} + u^{\sigma^2} e_{112}$$
$$+ u^\sigma u^{\sigma^2} e_{122} + u u^{\sigma^2} e_{212} + u u^\sigma e_{221} + \mathrm{N}_{k_1/k}(u) e_{222},$$
$$n(u) e_{112} = e_{112} + u e_{212} + u^\sigma e_{122} + u u^\sigma e_{222},$$
$$n(u) e_{121} = e_{121} + u e_{221} + u^{\sigma^2} e_{122} + u u^{\sigma^2} e_{222},$$
$$n(u) e_{122} = e_{122} + u e_{222},$$
$$n(u) e_{211} = e_{211} + u^\sigma e_{221} + u^{\sigma^2} e_{212} + u^\sigma u^{\sigma^2} e_{222},$$
$$n(u) e_{212} = e_{212} + u^\sigma e_{222},$$
$$n(u) e_{221} = e_{221} + u^{\sigma^2} e_{222},$$
$$n(u) e_{222} = e_{222}.$$

Therefore, if $n(u)x = \sum_{i,j,k} y_{ijk} e_{ijk}$,

$$(3.7) \qquad y_{111} = x_{111},$$
$$y_{211} = x_{211} + x_{111} u,$$
$$y_{122} = x_{122} + x_{111} u^\sigma u^{\sigma^2} + x_{211}^{\sigma^2} u^\sigma + x_{211}^\sigma u^{\sigma^2},$$
$$y_{222} = x_{222} + x_{111} \mathrm{N}_{k_1/k}(u) + \mathrm{Tr}_{k_1/k}(x_{211} u^\sigma u^{\sigma^2}) + \mathrm{Tr}_{k_1/k}(x_{122} u).$$

The element $\tau x$ is obtained by exchanging 1 and 2 in the indices of $e_{ijk}$'s (for example $\tau e_{122} = e_{211}$).

The relative invariant of $(G, V)$ can be constructed in the following manner. For $x = (x_{ijk}) \in \widetilde{V}$, we associate a $2 \times 2$ matrix with entries in the space of linear forms in two variables $v = (v_1, v_2)$ as

$$x \to M_x(v) = v_1 \begin{pmatrix} x_{111} & x_{121} \\ x_{211} & x_{221} \end{pmatrix} + v_2 \begin{pmatrix} x_{112} & x_{122} \\ x_{212} & x_{222} \end{pmatrix}.$$

Then $F_x(v) = \det M_x(v)$ is a quadratic form in $v = (v_1, v_2)$. Let $\Delta(x)$ be the discriminant of $F_x(v)$.

It was shown in [7] that

$$\Delta((g_1, g_2, g_3)x) = (\det g_1 \det g_2 \det g_3)^2 \Delta(x)$$

for $g_1, g_2, g_3 \in \mathrm{GL}(2)_{k_2}$ and $x \in \widetilde{V}$. So if we put $\chi(g) = \mathrm{N}_{k_1/k}(\det g)$ for $g \in \mathrm{GL}(2)_{k_1}$, $\chi$ is a $k$–rational character of $G$ and $\Delta(\phi(g)x) = \chi(g)^2 \Delta(x)$.

**Proposition (3.8)** $\Delta(x) \in k[V]$.

*Proof.* By an easy computation,

$$F_x(v) = (x_{111}x_{221} - x_{121}x_{211})v_1^2$$
$$+ (x_{111}x_{222} + x_{221}x_{112} - x_{121}x_{212} - x_{211}x_{122})v_1v_2$$
$$+ (x_{112}x_{222} - x_{122}x_{212})v_2^2.$$

So
$$\Delta(x) = (x_{111}x_{222} + x_{221}x_{112} - x_{121}x_{212} - x_{211}x_{122})^2$$
$$- 4(x_{111}x_{221} - x_{121}x_{211})(x_{112}x_{222} - x_{122}x_{212})$$
$$= x_{111}^2 x_{222}^2 + x_{221}^2 x_{112}^2 + x_{121}^2 x_{212}^2 + x_{211}^2 x_{122}^2$$
$$- 2x_{111}x_{222}x_{221}x_{112} - 2x_{111}x_{222}x_{121}x_{212}$$
$$- 2x_{111}x_{222}x_{211}x_{122} - 2x_{121}x_{221}x_{112}x_{212}$$
$$- 2x_{211}x_{221}x_{112}x_{122} - 2x_{121}x_{211}x_{122}x_{212}$$
$$+ 4x_{111}x_{221}x_{122}x_{212} + 4x_{121}x_{211}x_{112}x_{222}.$$

If $x \in V_k$,

$$\Delta(x) = x_{111}^2 x_{222}^2 + \mathrm{Tr}_{k_1/k}(x_{122}^2 x_{211}^2)$$
$$- 2x_{111}x_{222}\mathrm{Tr}_{k_1/k}(x_{122}x_{211}) - 2\mathrm{Tr}_{k_1/k}(x_{122}x_{122}^\sigma x_{211}x_{211}^\sigma)$$
$$+ 4x_{111}\mathrm{N}_{k_1/k}(x_{122}) + 4x_{222}\mathrm{N}_{k_1/k}(x_{211}).$$

This proves the proposition.

Q.E.D.

For the rest of this section, we consider an extra $\mathrm{GL}(1)_k$–factor because it is more natural number theoretically. So the group is $G = \mathrm{GL}(1)_k \times \mathrm{GL}(2)_{k_1}$ instead of $\mathrm{GL}(2)_{k_1}$. We can define an action of $G$ on $V$ by assuming that $t \in \mathrm{GL}(1)_k$ acts by the ordinary multiplication of $t$. Then $(G, V)$ is also a prehomogeneous vector space. Since the group is bigger, Condition (1.5) is still satisfied.

Let $w = e_{111} + e_{222}$. We identify $\tau$ with $(1, \tau)$. If $k_2 \subset k'$, the group $G^0_{wk'}$ is generated elements of the form

$$(t, a_2(t_{11}, t_{12}), a_2(t_{21}, t_{22}), a_2(t_{31}, t_{32})),$$

and $G_{wk'}$ is generated by $G^0_{wk'}$ and $\tau$. Therefore we get the following proposition.

**Proposition (3.9)** (1) *As an algebraic group over $k$,*

$$G_w^0 = \left\{(t_1, a_2(t_{21}, t_{22})) \,\middle|\, \begin{array}{c} t_1 \in k^\times,\ t_{21}, t_{22} \in k_1^\times, \\ t_1 \mathrm{N}_{k_1/k}(t_{21}) = t_1 \mathrm{N}_{k_1/k}(t_{22}) = 1 \end{array}\right\}.$$

(2) *There is a split exact sequence*

$$1 \to G_w^0 \to G_w \to \mathfrak{S}_2 \to 1$$

*where the action of the Galois group on $\mathfrak{S}_2$ is trivial.*

From (1.8) we have a surjective map $\alpha_V : G_k \backslash V_k^{\mathrm{ss}} \to \mathfrak{Ex}_2$.

Let $k(\alpha)/k$ be a separable quadratic extension. Then $k_1 \otimes k(\alpha) \cong k_1(\alpha)$. If $k_2 \not\supset k(\alpha)$, we extend $\sigma$ to $k_1(\alpha)$ so that it is trivial on $k(\alpha)$. If $k_2/k$ is an $\mathfrak{S}_3$–extension containing $k(\alpha)$, $k(\alpha)$ is the unique quadratic extension of $k$ contained in $k_2$. So $\sigma$ is trivial on $k(\alpha)$. In both cases, we can regard the action of $G_{k(\alpha)} = \mathrm{GL}(1)_{k(\alpha)} \times \mathrm{GL}(2)_{k_1(\alpha)}$ on $V_{k(\alpha)}$ as given by

$$(t, g) \cdot x_1 \otimes x_2 \otimes x_3 = tgx_1 \otimes g^\sigma x_2 \otimes g^{\sigma^2} x_3$$

for $t \in k(\alpha)$, $g \in \mathrm{GL}(2)_{k_1(\alpha)}$.

Let $\nu$ be the non-trivial element of $\mathrm{Gal}(k(\alpha)/k)$. We extend $\nu$ to $\mathrm{Gal}(k_2(\alpha)/k)$ so that it is trivial on $k_1$. Note that this is possible even if $\alpha \in k_2$, because if $k_1$ is generated by $\alpha_1 \in k_1$ over $k$ and $\alpha_2, \alpha_3 \in k_2$ are conjugate elements, we can choose $\nu$ to be the transposition (23). By this extension, $\nu$ induces an involution on

$$G_{k(\alpha)} \cong \mathrm{GL}(1)_{k(\alpha)} \times \mathrm{GL}(2)_{k_1(\alpha)},$$

where the action of $\nu$ on the right hand side is defined by the extension of $\nu$ to $\mathrm{Gal}(k_2(\alpha)/k)$, and

$$\{g \in G_{k(\alpha)} \mid g^\nu = g\} = G_k.$$

Let $f(v) = v_1^2 + a_1 v_1 v_2 + a_2 v_2^2 \in k[v]$ be a polynomial such that the roots $\alpha_1, \alpha_2$ of $f$ generate the field $k(\alpha)$. We define

$$(3.10) \quad g_\alpha = \left(1, \begin{pmatrix} 1 & 1 \\ -\alpha_1 & -\alpha_2 \end{pmatrix}\right),$$

$$w_\alpha = 2e_{111} + a_1(e_{211} + e_{121} + e_{112})$$
$$+ (a_1^2 - 2a_2)(e_{122} + e_{212} + e_{221}) + (a_1^3 - 3a_1 a_2)e_{222}.$$

Then $g_\alpha w = w_\alpha$. Since $g_\alpha^\nu = g_\alpha \tau$ for the non-trivial element $\nu$ of $\mathrm{Gal}(k(\alpha)/k)$, $k(\alpha)$ corresponds to the cohomology class $\{g_\alpha^\eta g_\alpha^{-1}\} \in \mathrm{H}^1(k, \mathfrak{S}_2)$. Therefore, $k(w_\alpha) = k(\alpha)$.

For the rest of this section, we determine $G_x^0$ for all $x \in V_k^{\mathrm{ss}}$.

By the definition of $w_\alpha$,

$$G_{w_\alpha k^{\mathrm{sep}}} = g_\alpha G_{w k^{\mathrm{sep}}} g_\alpha^{-1}, \quad G_{w_\alpha k^{\mathrm{sep}}}^0 = g_\alpha G_{w k^{\mathrm{sep}}}^0 g_\alpha^{-1}.$$

Note that
$$(1, a_2(1, -1)) = g_\alpha \tau g_\alpha^{-1} \in G_{w_\alpha k}.$$
Therefore, we have a split exact sequence

(3.11) $$1 \to G^0_{w_\alpha} \to G_{w_\alpha} \to \mathfrak{S}_2 \to 1,$$

where the Galois group acts on $\mathfrak{S}_2$ trivially.

We get the following proposition by the same argument as in §2.

**Proposition (3.12)** *As an algebraic group over $k$,*

$$G^0_{w_\alpha} = \left\{ g_\alpha(t_1, a_2(t_{21}, t_{21}^\nu))g_\alpha^{-1} \,\middle|\, \begin{array}{c} t_1 \in k^\times,\ t_{21} \in k_1(\alpha)^\times, \\ t_1 N_{k_1(\alpha)/k(\alpha)}(t_{21}) = t_1 N_{k_1(\alpha)/k(\alpha)}(t_{21}^\nu) = 1 \end{array} \right\}$$
$$\cong \{t_{21} \in \mathrm{GL}(1)_{k_1(\alpha)} \mid N_{k_1(\alpha)/k(\alpha)}(t_{21}) \in \mathrm{GL}(1)_k\}.$$

**Theorem (3.13)** *There is a bijection*

$$\mathfrak{S}_2 \setminus (k^\times / N_{k_1/k}(k_1^\times))^2 / \{(t, t) \mid t \in k^\times\} \cong \alpha_V^{-1}(k),$$

*where the action of $\mathfrak{S}_2$ is given by permutations. Moreover $\beta = (\beta_1, \beta_2) \in (k^\times)^2$ corresponds to the orbit of $\beta_1 e_{111} + \beta_2 e_{222}$.*
*(2) If $k(\alpha)/k$ is a quadratic extension, there is a bijection*

$$\mathfrak{S}_2 \setminus (k(\alpha)^\times / k^\times N_{k_1(\alpha)/k(\alpha)}(k_1(\alpha)^\times)) \cong \alpha_V^{-1}(k(\alpha)),$$

*where the action of the non-trivial element $\nu \in \mathfrak{S}_2 \cong \mathrm{Gal}(k(\alpha)/k)$ is given by the usual Galois group action. Moreover, $\beta \in k(\alpha)^\times$ corresponds to the orbit of $g_\alpha(\beta e_{111} + \beta^\nu e_{222})$.*

*Proof.* We have an exact sequence

$$1 \to G^0_w \to \mathrm{GL}(1)_{k_1} \times \mathrm{GL}(1)_{k_1} \to \mathrm{GL}(1)_k \to 1,$$

where the last map is given by $(t_1, t_2) \to N_{k_1/k}(t_1 t_2^{-1})$. So we have an exact sequence

$$k_1^\times \times k_1^\times \to k^\times \to \mathrm{H}^1(k, G^0_w) \to 1.$$

This implies

$$\mathrm{H}^1(k, G^0_w) \cong k^\times / N_{k_1/k}(k_1^\times) \cong (k^\times / N_{k_1/k}(k_1^\times))^2 / \{(t, t) \mid t \in k^\times\}.$$

We will calculate the orbit corresponding to an element of $\mathrm{H}^1(k, G^0_w)$ realized in this way. This amounts to making explicit the boundary map in the above sequence. Let $\beta = (\beta_1, \beta_2) \in (k^\times)^2$. We choose a a large enough finite Galois extension $k'/k$ and $\beta' = (\beta'_1, \beta'_2) \in k_1 \otimes k' \times k_1 \otimes k'$ such that $N_{k_1/k}(\beta'_i) = \beta_i$ for $i = 1, 2$. Here we are considering $N_{k_1/k}$ as a $k$–morphism between the $k$–varieties $\mathrm{GL}(1)_{k_1}$ and $\mathrm{GL}(1)_k$. Then the corresponding element in $\mathrm{H}^1(k, G_w^0)$ is given by $\{\beta'^{-1}\beta'^\eta\}_{\eta\in\mathrm{Gal}(k'/k)}$. By Proposition (3.9), this corresponds to the element

(3.14) $$(1, a_2(\beta'_1, \beta'_2))^{-1}(1, a_2(\beta'_1, \beta'_2))^\eta.$$

Let $g_\beta = (1, a_2(\beta'_1, \beta'_2)) \in G_{k'}$. Then in $\mathrm{H}^1(k'/k, G_{k'})$, the element (3.14) is trivial and becomes $\delta g_\beta$. Therefore, $\beta$ corresponds to the orbit of

$$\begin{aligned} g_\beta w &= \mathrm{N}_{k_1/k}(\beta'_1) e_{111} + \mathrm{N}_{k_1/k}(\beta'_2) e_{222} \\ &= \beta_1 e_{111} + \beta_2 e_{222}. \end{aligned}$$

Clearly $\tau$ exchanges $\beta_1, \beta_2$. This proves (1).

Consider (2). Let

$$A = \{t \in \mathrm{GL}(1)_{k_1(\alpha)} \mid \mathrm{N}_{k_1(\alpha)/k(\alpha)}(t) = 1\}.$$

Consider the following two exact sequences of abelian groups

$$\begin{array}{ccccccccc} 1 & \to & A & \to & G_{w_\alpha}^0 & \to & \mathrm{GL}(1)_k & \to & 1 \\ & & \| & & \downarrow & & \downarrow & & \\ 1 & \to & A & \to & \mathrm{GL}(1)_{k_1(\alpha)} & \to & \mathrm{GL}(1)_{k(\alpha)} & \to & 1 \end{array}$$

where the last maps are given by $\mathrm{N}_{k_1(\alpha)/k(\alpha)}$.

From the above exact sequences, we have the following long exact sequences.

$$\begin{array}{ccccccccc} \cdots & \to & k^\times & \to & \mathrm{H}^1(k, A) & \to & \mathrm{H}^1(k, G_{w_\alpha}^0) & \to & 1 \\ & & \downarrow & & \| & & & & \\ k_1(\alpha)^\times & \to & k(\alpha)^\times & \to & \mathrm{H}^1(k, A) & \to & & & 1 \end{array}$$

So,
$$\mathrm{H}^1(k, G_{w_\alpha}^0) \cong \mathrm{H}^1(k, A)/k^\times,$$
$$\mathrm{H}^1(k, A) \cong k(\alpha)^\times / \mathrm{N}_{k_1(\alpha)/k(\alpha)}(k_1(\alpha)^\times).$$

Therefore,
$$\mathrm{H}^1(k, G_{w_\alpha}^0) \cong k(\alpha)^\times / k^\times \mathrm{N}_{k_1(\alpha)/k(\alpha)}(k_1(\alpha)^\times).$$

Let $\beta \in k(\alpha)^\times$. We consider the image of $\beta$ by the boundary map $k(\alpha)^\times \to \mathrm{H}^1(k, G_{w_\alpha}^0)$. We choose a large enough finite Galois extension $k'/k$ and $\beta' \in k_1(\alpha) \otimes k'$ such that $\beta = \mathrm{N}_{k_1(\alpha)/k(\alpha)}(\beta')$. Here we are considering $\mathrm{N}_{k_1(\alpha)/k(\alpha)}$ as a $k$-morphism between the $k$-varieties $\mathrm{GL}(1)_{k_1(\alpha)}$ and $\mathrm{GL}(1)_{k(\alpha)}$. Then the corresponding element in $\mathrm{H}^1(k, G_{w_\alpha}^0)$ is given by $\{\beta'^{-1}\beta'^\eta\}_{\eta\in\mathrm{Gal}(k'/k)}$. Note that since we are regarding $\beta$ as a $k$-rational point of $\mathrm{GL}(1)_{k(\alpha)}$, $\beta^\eta = \beta$ for all $\eta$. So $\mathrm{N}_{k_1(\alpha)/k(\alpha)}(\beta'^{-1}\beta'^\eta) = \beta^{-1}\beta^\eta = 1$. By Proposition (3.12), $\beta'^{-1}\beta'^\eta$ corresponds to

(3.15) $$g_\alpha(1, a_2(\beta'^{-1}\beta'^\eta, (\beta'^{-1}\beta'^\eta)^\nu))g_\alpha^{-1} \in G_{wk'}^0.$$

Here $\nu \in \mathrm{Gal}(k(\alpha)/k)$ acts on $k(\alpha) \otimes k'$ by the first factor and we are regarding $\nu$ as a $k$-automorphism of the $k$-varieties $\mathrm{GL}(1)_{k_1(\alpha)}$ and $\mathrm{GL}(1)_{k(\alpha)}$. But $\eta$ acts on $k(\alpha) \otimes k'$ by the second factor. Therefore $\nu$ and $\eta$ commute, even though they need not in $\mathrm{Gal}(k^{\mathrm{sep}}/k)$, and we are slightly abusing notation. Similarly we are regarding $g_\alpha$ as a point having coordinates in $k(\alpha) \otimes k'$. Therefore $g_\alpha$ is also fixed by $\eta$ for all $\eta$. So we can write the element (3.15) as

$$(g_\alpha(1, a_2(\beta', \beta'^\nu))g_\alpha^{-1})^{-1}(g_\alpha(1, a_2(\beta', \beta'^\nu))g_\alpha^{-1})^\eta.$$

Let $g_\beta = g_\alpha(1, a_2(\beta', \beta'^\nu))g_\alpha^{-1} \in G_{k'}$. Then in $\mathrm{H}^1(k, G_{k'})$, the above 1–cocycle is trivial and becomes $\delta g_\beta$. Therefore, $\beta$ corresponds to the orbit of

$$\begin{aligned}
g_\beta w_\alpha &= g_\alpha(1, a_2(\beta', \beta'^\nu))g_\alpha^{-1} g_\alpha w \\
&= g_\alpha(1, a_2(\beta', \beta'^\nu))w \\
&= g_\alpha(\mathrm{N}_{k_1(\alpha)/k(\alpha)}(\beta')e_{111} + \mathrm{N}_{k_1(\alpha)/k(\alpha)}(\beta'^\nu)e_{222}) \\
&= g_\alpha(\beta e_{111} + \beta^\nu e_{222}).
\end{aligned}$$

Clearly, $g_\alpha \tau g_\alpha^{-1}$ maps $\beta$ to $\beta^\nu$. This proves the theorem. Note if we can write the element (3.15) as

$$(g_\alpha(\beta^{-1}, a_2(\beta', \beta'^\nu))g_\alpha^{-1})^{-1}(g_\alpha(\beta^{-1}, a_2(\beta', \beta'^\nu))g_\alpha^{-1})^\eta,$$

it determines a trivial class in $\mathrm{H}^1(k, G^0_{w_\alpha})$. But this is not necessarily possible because $\beta$ is not necessarily a $k$–rational point of $\mathrm{GL}(1)_k$.

Q.E.D.

**Corollary (3.16)** (1) If $x, y \in V_k^{\mathrm{ss}}$ and $k(x) = k(y)$ then $G_x^0 = G_y^0$.
(2) If the characteristic of $k$ is not two, $k(x) = k(\Delta(x)^{\frac{1}{2}})$ for all $x \in V_k^{\mathrm{ss}}$.

*Proof.* The statement (1) is clear from the description of each orbit.
Consider (2). Suppose that the characteristic of $k$ is not two.
For $\beta$ in Theorem (3.13)(1),

$$\Delta(\beta_1 e_{111} + \beta_2 e_{222}) = (\beta_1 \beta_2)^2,$$

and for this element, the field extension is trivial.
For $\beta$ in Theorem (3.13)(2),

$$\begin{aligned}
\Delta(g_\beta w_\alpha) &= \Delta(g_\alpha(\beta e_{111} + \beta^\nu e_{222}) \\
&= \chi(g_\alpha)^2 \Delta(\beta e_{111} + \beta^\nu e_{222}) \\
&= \chi(g_\alpha)^2 \mathrm{N}_{k(\alpha)/k}(\beta)^2 \\
&= \Delta(w_\alpha) \mathrm{N}_{k(\alpha)/k}(\beta)^2 = (\alpha_1 - \alpha_2)^6 \mathrm{N}_{k(\alpha)/k}(\beta)^2.
\end{aligned}$$

Since the characteristic of $k$ is not two,

$$k(g_\beta w_\alpha) = k(w_\alpha) = k(\alpha) = k(\alpha_1 - \alpha_2) = k(\Delta(g_\beta w_\alpha)^{\frac{1}{2}}).$$

This proves (2).

Q.E.D.

## §4. The non-split $E_6$ case

In this section, we consider the space of pairs of ternary Hermitian forms. We fix a separable quadratic extension $k_1 = k(\alpha_0)$ of $k$. The non-trivial element of $\mathrm{Gal}(k_1/k)$ is denoted by $\sigma$. Let $G = \mathrm{GL}(3)_{k_1} \times \mathrm{GL}(2)_k$ considered as an algebraic group over $k$. Let $W$ be the space of ternary Hermitian forms, and $V = W \otimes k^2$. The definitions of the action of $G$, $F_x(v)$, $\mathrm{Zero}(x)$, $\Delta(x)$, $k(x)$ for $x \in V_k^{\mathrm{ss}}$, etc., are similar to those in §2.

By a similar argument as in §2, we can show that $(G, V)$ becomes the $E_6$ case in [7]. Therefore, $(G, V)$ is a prehomogeneous vector space. It is possible to check by linear algebra that the $E_6$ case in [7] satisfies Condition (1.5). Since $k_1/k$ is a separable extension, we get the following proposition.

**Proposition (4.1)** *The prehomogeneous vector space $(G, V)$ satisfies* Condition (1.5).

Let

$$w = v_1 \begin{pmatrix} 1 & & \\ & -1 & \\ & & 0 \end{pmatrix} + v_2 \begin{pmatrix} 0 & & \\ & -1 & \\ & & 1 \end{pmatrix},$$

$$\tau_1 = \left( \begin{pmatrix} 0 & 1 & 0 \\ 1 & 0 & 0 \\ 0 & 0 & 1 \end{pmatrix}, \begin{pmatrix} -1 & 0 \\ -1 & 1 \end{pmatrix} \right),$$

$$\tau_2 = \left( \begin{pmatrix} 0 & 0 & 1 \\ 0 & 1 & 0 \\ 1 & 0 & 0 \end{pmatrix}, \begin{pmatrix} 0 & 1 \\ 1 & 0 \end{pmatrix} \right).$$

Then $\tau_1, \tau_2 \in G_{wk}$.

It is easy to see that $\mathrm{Zero}(w) = \{(1, -1), (0, 1), (1, 0)\}$ and $\tau_1, \tau_2$ correspond to the permutations (12) and (13). By the same argument as in §2 the sequence (1.9) is exact for all $x \in V_k^{\mathrm{ss}}$. Since $\tau_1, \tau_2$ are rational elements, $\mathrm{Aut}(\mathrm{Zero}(w)) \cong \mathfrak{S}_3$. Therefore, we have a split exact sequence

(4.2) $$1 \to G_w^0 \to G_w \to \mathfrak{S}_3 \to 1,$$

where the action of the Galois group on $\mathfrak{S}_3$ is trivial. So we have a surjective map $\alpha_V : G_k \backslash V_k^{\mathrm{ss}} \to \mathfrak{E}\mathfrak{x}_3 \cong \mathrm{H}^1(k, \mathfrak{S}_3)$.

By considering $G_{w k^{\mathrm{sep}}}^0$, we get the following proposition.

**Proposition (4.3)** *As an algebraic group over $k$*,

$$G_w^0 = \left\{ (a_3(t_{11}, t_{12}, t_{13}), t_2 I_2) \,\middle|\, \begin{array}{c} t_{11}, t_{12}, t_{13} \in k_1^\times,\ t_2 \in k^\times, \\ t_2 \mathrm{N}_{k_1/k}(t_{11}) = t_2 \mathrm{N}_{k_1/k}(t_{12}) = t_2 \mathrm{N}_{k_1/k}(t_{13}) = 1 \end{array} \right\}.$$

Let $f(v) = v_1^3 + a_1 v_1^2 v_2 + a_2 v_1 v_2^2 + a_3 v_2^3 \in k[v]$ be a cubic polynomial without a multiple factor and $\alpha_1, \alpha_2, \alpha_3$ the roots of $f(v)$.

We define

(4.4)
$$D_\alpha = (\alpha_1 - \alpha_2)(\alpha_1 - \alpha_3)(\alpha_2 - \alpha_3),$$
$$Q_\alpha = \frac{1}{D_\alpha}\begin{pmatrix} -(\alpha_2 - \alpha_3) & \alpha_2 - \alpha_1 \\ \alpha_1(\alpha_2 - \alpha_3) & -\alpha_3(\alpha_2 - \alpha_1) \end{pmatrix},$$
$$P_\alpha = \begin{pmatrix} 1 & 1 & 1 \\ \alpha_1 & \alpha_2 & \alpha_3 \\ \alpha_1^2 & \alpha_2^2 & \alpha_3^2 \end{pmatrix},$$
$$w_\alpha = v_1 \begin{pmatrix} 0 & 0 & -1 \\ 0 & -1 & a_1 \\ -1 & a_1 & -a_1^2 + a_2 \end{pmatrix}$$
$$+ v_2 \begin{pmatrix} 0 & 1 & -a_1 \\ 1 & -a_1 & a_1^2 - a_2 \\ -a_1 & a_1^2 - a_2 & -a_1^3 + 2a_1 a_2 - a_3 \end{pmatrix}.$$

By easy computations,

$$\det P_\alpha = D_\alpha, \quad \det Q_\alpha = \frac{1}{D_\alpha},$$
$$\det w_\alpha = f(v).$$

Let $k(\alpha) = k(\alpha_1, \alpha_2, \alpha_3)$ and $k_1(\alpha) = k_1(\alpha_1, \alpha_2, \alpha_3)$. Then $k(\alpha)$ and $k_1(\alpha)$ are Galois extensions of $k$. By a similar argument as in §2, $k(w_\alpha) = k(\alpha)$.

If $k_1 \not\subset k(\alpha)$,
$$G_{k(\alpha)} = \mathrm{GL}(3)_{k_1(\alpha)} \times \mathrm{GL}(2)_{k(\alpha)},$$
$$V_{k(\alpha)} = W \otimes k(\alpha)^2.$$

We extend $\sigma$ to $\mathrm{Gal}(k_1(\alpha)/k)$ so that it is trivial on $k(\alpha)$. If $\nu \in \mathrm{Gal}(k(\alpha)/k)$ we extend $\nu$ to $\mathrm{Gal}(k_1(\alpha)/k)$ so that it is trivial on $k_1$. Then $\nu$ acts $G_{k(\alpha)}$ and

$$G_k = \{g \in G_{k(\alpha)} \mid g^\nu = g \text{ for all } \nu \in \mathrm{Gal}(k(\alpha)/k)\}.$$

If $k_1 \subset k(\alpha)$,

$$G_{k(\alpha)} = \mathrm{GL}(3)_{k(\alpha)} \times \mathrm{GL}(3)_{k(\alpha)} \times \mathrm{GL}(2)_{k(\alpha)},$$
$$V_{k(\alpha)} = \mathrm{M}(2,2) \otimes k(\alpha)^2.$$

If $k(\alpha) = k_1$, the action of $\sigma$ on $G_{k(\alpha)}$ is given by

$$(g_1, g_2, g_3)^\sigma = (g_2^\sigma, g_1^\sigma, g_3^\sigma).$$

Also the action of $G_{k(\alpha)}$ on $V_{k(\alpha)}$ is given by

$$(g_1, g_2, g_3) M(v) = g_1 M(v g_3)^t g_2.$$

If $k(\alpha)/k$ is an $\mathfrak{S}_3$-extension containing $k_1$, we extend $\sigma$ to $\mathrm{Gal}(k(\alpha)/k)$ so that it is trivial on $k(\alpha_1)$. In other words we are regarding $\sigma$ as (23). The action of $\nu \in \mathrm{Gal}(k(\alpha)/k)$ is given by

$$(g_1, g_2, g_3)^\nu = \begin{cases} (g_1^\nu, g_2^\nu, g_3^\nu) & \nu \text{ is trivial on } k_1, \\ (g_2^\nu, g_1^\nu, g_3^\nu) & \nu \text{ is not trivial on } k_1. \end{cases}$$

Also

$$G_k = \{g \in G_{k(\alpha)} \mid g^\nu = g \text{ for all } \nu \in \mathrm{Gal}(k(\alpha)/k)\}.$$

We define $g_\alpha \in G_{k^{\mathrm{sep}}}$ by

(4.5)
$$g_\alpha = \begin{cases} (P_\alpha, Q_\alpha) & k_1 \not\subset k(\alpha), \\ (P_\alpha, P_\alpha, Q_\alpha) & k_1 \subset k(\alpha). \end{cases}$$

**Lemma (4.6)** $w_\alpha = g_\alpha w$.

*Proof.* Let

$$A_i(\alpha) = \alpha_1^i(\alpha_2 - \alpha_3) + \alpha_2^i(\alpha_3 - \alpha_1) + \alpha_3^i(\alpha_1 - \alpha_2)$$

for $i = 2, 3, 4, 5$. We define

$$W_1(\alpha) = \begin{pmatrix} -(\alpha_2 - \alpha_3) & & \\ & -(\alpha_3 - \alpha_1) & \\ & & -(\alpha_1 - \alpha_2) \end{pmatrix},$$

$$W_2(\alpha) = \begin{pmatrix} -\alpha_1(\alpha_2 - \alpha_3) & & \\ & -\alpha_2(\alpha_3 - \alpha_1) & \\ & & -\alpha_3(\alpha_1 - \alpha_2) \end{pmatrix},$$

$$W_3(\alpha) = \begin{pmatrix} 0 & 0 & -A_2(\alpha) \\ 0 & -A_2(\alpha) & -A_3(\alpha) \\ -A_2(\alpha) & -A_3(\alpha) & -A_4(\alpha) \end{pmatrix},$$

$$W_4(\alpha) = \begin{pmatrix} 0 & A_2(\alpha) & A_3(\alpha) \\ A_2(\alpha) & A_3(\alpha) & A_4(\alpha) \\ A_3(\alpha) & A_4(\alpha) & A_5(\alpha) \end{pmatrix}.$$

Then

$$g_\alpha w = \frac{1}{D_\alpha} P_\alpha (v_1 W_1(\alpha) + v_2 W_2(\alpha))^t P_\alpha$$

$$= \frac{1}{D_\alpha}(v_1 W_3(\alpha) + v_2 W_4(\alpha)).$$

Now the lemma follows from the relations

$$A_2(\alpha) = D_\alpha,$$
$$A_3(\alpha) = -D_\alpha a_1,$$
$$A_4(\alpha) = D_\alpha(a_1^2 - a_2),$$
$$A_5(\alpha) = D_\alpha(-a_1^3 + 2a_1 a_2 - a_3).$$

Q.E.D.

Let

$$\text{(4.7)} \quad t = \begin{cases} (a_3(t_{11},t_{12},t_{13}), t_2 I_2) & k_1 \not\subset k(\alpha), \\ (a_3(t_{11},t_{12},t_{13}), a_3(t_{21},t_{22},t_{23}), t_3 I_2) & k_1 \subset k(\alpha). \end{cases}$$

In the following proposition if $k(\alpha)/k$ is a quadratic extension different from $k_1$ then $\nu \in \operatorname{Gal}(k(\alpha)/k) \cong \operatorname{Gal}(k_1(\alpha)/k_1)$ is the non-trivial element and if $k(\alpha_1)/k$ is a cubic extension, $\nu \in \operatorname{Gal}(k_1(\alpha)/k)$ is the element such that $\nu(\alpha_1) = \alpha_2$, $\nu(\alpha_2) = \alpha_3$, and $\nu(s) = s$ for all $s \in k_1$. Also if $k(\alpha)/k$ is a quadratic extension, we choose $\alpha_1, \alpha_2$ conjugate over $k$ and $\alpha_3 = 0$.

**Proposition (4.8)** (1) *If $k(\alpha)$ is a quadratic extension of $k$ different from $k_1$,*

$$G^0_{w_\alpha} = \left\{ g_\alpha t g_\alpha^{-1} \,\middle|\, \begin{array}{c} t_{11}, t_{12} \in k_1(\alpha)^\times, t_{13} \in k_1^\times, t_2 \in k^\times, t_{12} = t_{11}^\nu \\ t_2 \mathrm{N}_{k_1(\alpha)/k(\alpha)}(t_{11}) = t_2 \mathrm{N}_{k_1(\alpha)/k(\alpha)}(t_{12}) = t_2 \mathrm{N}_{k_1/k}(t_{13}) = 1 \end{array} \right\}$$
$$\cong \operatorname{Ker}(\operatorname{GL}(1)_{k_1(\alpha)} \times \operatorname{GL}(1)_{k_1} \to \operatorname{GL}(1)_{k(\alpha)}),$$

*where $\operatorname{GL}(1)_{k_1(\alpha)} \times \operatorname{GL}(1)_{k_1} \to \operatorname{GL}(1)_{k(\alpha)}$ is given by*

$$(t_{11}, t_{13}) \to \mathrm{N}_{k_1(\alpha)/k(\alpha)}(t_{11} t_{13}^{-1}).$$

(2) *If $k(\alpha) = k_1$,*

$$G^0_{w_\alpha} = \left\{ g_\alpha t g_\alpha^{-1} \,\middle|\, \begin{array}{c} t_{ij} \in k_1^\times \text{ for } i=1,2, j=1,2,3, t_3 \in k^\times, \\ t_{22} = t_{11}^\sigma, t_{21} = t_{12}^\sigma, t_{23} = t_{13}^\sigma, t_3 t_{11} t_{12}^\sigma = t_3 \mathrm{N}_{k_1/k}(t_{13}) = 1 \end{array} \right\}$$
$$\cong \operatorname{GL}(1)_{k_1} \times \operatorname{GL}(1)_{k_1}.$$

(3) *If $k(\alpha)$ is either a cyclic cubic extension or an $\mathfrak{S}_3$-extension of $k$ not containing $k_1$,*

$$G^0_{w_\alpha} = \left\{ g_\alpha t g_\alpha^{-1} \,\middle|\, \begin{array}{c} t_{1j} \in k_1(\alpha_j)^\times \text{ for } j=1,2,3, t_2 \in k^\times, \\ t_{12} = t_{11}^\nu, t_{13} = t_{11}^{\nu^2}, t_2 \mathrm{N}_{k_1(\alpha_j)/k(\alpha_j)}(t_{1j}) = 1 \text{ for } j=1,2,3 \end{array} \right\}$$
$$\cong \{t_{11} \in \operatorname{GL}(1)_{k_1(\alpha_1)} \mid \mathrm{N}_{k_1(\alpha_1)/k(\alpha_1)}(t_{11}) \in \operatorname{GL}(1)_k \}.$$

(4) *If $k(\alpha)$ is an $\mathfrak{S}_3$-extension of $k$ containing $k_1$,*

$$G^0_{w_\alpha} = \left\{ g_\alpha t g_\alpha^{-1} \,\middle|\, \begin{array}{c} t_{ij} \in k(\alpha)^\times \text{ for } i=1,2, j=1,2,3, t_3 \in k^\times, \\ t_{12} = t_{11}^\nu, t_{13} = t_{11}^{\nu^2}, t_{21} = t_{11}^\sigma, t_{22} = t_{11}^{\sigma\nu}, t_{23} = t_{11}^{\sigma\nu^2} \\ t_3 \mathrm{N}_{k(\alpha)/k(\alpha_1)}(t_{11}) = t_3 \mathrm{N}_{k(\alpha)/k(\alpha_1)}(t_{21}) = 1 \end{array} \right\}$$
$$\cong \{t_{11} \in \operatorname{GL}(1)_{k(\alpha)} \mid \mathrm{N}_{k(\alpha)/k(\alpha_1)}(t_{11}) \in \operatorname{GL}(1)_k \}.$$

*Proof.* In the following proof, we only consider the set $G^0_{w_\alpha k}$ of $k$-rational points, but the argument can easily be generalized to $G^0_{w_\alpha R}$ for any $k$-algebra $R$ as in Proposition (2.10). Therefore, we are proving isomorphisms of algebraic groups over $k$.

Consider (1). It is easy to see that $g_\alpha^\nu = g_\alpha \tau_1$. So

$$(g_\alpha t g_\alpha^{-1})^\nu = g_\alpha \tau_1 t^\nu \tau_1 g_\alpha^{-1}$$
$$= g_\alpha (a_3(t_{12}^\nu, t_{11}^\nu, t_{13}^\nu), t_2^\nu I_2) g_\alpha^{-1}.$$

Therefore, $(g_\alpha t g_\alpha^{-1})^\nu = g_\alpha t g_\alpha^{-1}$ if and only if $t_{12} = t_{11}^\nu$, $t_{13} \in k_1^\times$, $t_2 \in k^\times$. The rest of the condition is obvious and this proves (1).

Consider (2). Note that in $G_{k(\alpha)}$, $\tau_1$ corresponds to the element

$$\left( \begin{pmatrix} 0 & 1 & 0 \\ 1 & 0 & 0 \\ 0 & 0 & 1 \end{pmatrix}, \begin{pmatrix} 0 & 1 & 0 \\ 1 & 0 & 0 \\ 0 & 0 & 1 \end{pmatrix}, \begin{pmatrix} -1 & 0 \\ -1 & 1 \end{pmatrix} \right).$$

Since $g_\alpha^\sigma = g_\alpha \tau_1$,

$$(g_\alpha t g_\alpha^{-1})^\sigma = g_\alpha \tau_1 t^\sigma \tau_1 g_\alpha^{-1}$$
$$= g_\alpha (a_3(t_{22}^\sigma, t_{21}^\sigma, t_{23}^\sigma), a_3(t_{12}^\sigma, t_{11}^\sigma, t_{13}^\sigma), t_3^\sigma I_2) g_\alpha^{-1}.$$

Therefore, $(g_\alpha t g_\alpha^{-1})^\sigma = g_\alpha t g_\alpha^{-1}$ if and only if $t_{21} = t_{12}^\sigma$, $t_{22} = t_{11}^\sigma$, $t_{23} = t_{13}^\sigma$, $t_3 \in k^\times$.

The rest of the condition is obvious, and this proves (2).

Consider (3). It is easy to see that $g_\alpha^\nu = g_\alpha \tau_2 \tau_1$. So

$$(g_\alpha t g_\alpha^{-1})^\nu = g_\alpha \tau_2 \tau_1 t^\nu \tau_1 \tau_2 g_\alpha^{-1}$$
$$= g_\alpha (a_3(t_{13}^\nu, t_{11}^\nu, t_{12}^\nu), t_2^\nu I_2) g_\alpha^{-1}.$$

Therefore, , $(g_\alpha t g_\alpha^{-1})^\nu = g_\alpha t g_\alpha^{-1}$ if and only if $t_{12} = t_{11}^\nu$, $t_{13} = t_{12}^\nu = t_{11}^{\nu^2}$, $t_2^\nu = t_2$.

If $k(\alpha)/k$ is a cyclic cubic extension, $k(\alpha) = k(\alpha_1) = k(\alpha_2) = k(\alpha_3)$. If $k(\alpha)/k$ is an $\mathfrak{S}_3$-extension, let $\nu' \in \mathrm{Gal}(k(\alpha)/k)$ be the element which corresponds to the transposition (12). We extend $\nu'$ to $\mathrm{Gal}(k_1(\alpha)/k)$ so that it is trivial on $k_1$. Since $g_\alpha^{\nu'} = g_\alpha \tau_1$,

$$(g_\alpha t g_\alpha^{-1})^{\nu'} = g_\alpha \tau_1 t^{\nu'} \tau_1 g_\alpha^{-1}$$
$$= g_\alpha (a_3(t_{12}^{\nu'}, t_{11}^{\nu'}, t_{13}^{\nu'}), t_2^{\nu'} I_2) g_\alpha^{-1}.$$

Therefore, $(g_\alpha t g_\alpha^{-1})^{\nu'} = g_\alpha t g_\alpha^{-1}$ if and only if $t_{12} = t_{11}^{\nu'}$, $t_{13} \in k_1(\alpha_3)$, $t_2^{\nu'} = t_2$. Similarly, we can prove that $t_{11} \in k_1(\alpha_1)$, $t_{12} \in k_1(\alpha_2)$.

The rest of the condition is obvious. Since $\mathrm{Gal}(k(\alpha)/k)$ is generated by the permutations (13) and (123), this proves (3).

Consider (4). Note that since $\nu$ is trivial on $k_1$, $(g_1, g_2, g_3)^\nu = (g_1^\nu, g_2^\nu, g_3^\nu)$ for $(g_1, g_2, g_3) \in G_{k(\alpha)}$. So, as in (3), $(g_\alpha t g_\alpha^{-1})^\nu = g_\alpha t g_\alpha^{-1}$ if and only if

$$t_{12} = t_{11}^\nu, \ t_{13} = t_{12}^\nu, \ t_{22} = t_{21}^\nu, \ t_{23} = t_{22}^\nu, \ t_3^\nu = t_3.$$

Note that we extended $\sigma$ to $\mathrm{Gal}(k(\alpha)/k)$ so that it corresponds to the transposition (23). Since $g_\alpha^\sigma = g_\alpha \tau_1 \tau_2 \tau_1$,

$$(g_\alpha t g_\alpha^{-1})^\sigma = g_\alpha \tau_1 \tau_2 \tau_1 t^\sigma \tau_1 \tau_2 \tau_1 g_\alpha^{-1}$$
$$= g_\alpha (a_3(t_{21}^\sigma, t_{23}^\sigma, t_{22}^\sigma), a_3(t_{11}^\sigma, t_{13}^\sigma, t_{12}^\sigma), t_3^\sigma I_2) g_\alpha^{-1}.$$

Therefore, , $(g_\alpha t g_\alpha^{-1})^\sigma = g_\alpha t g_\alpha^{-1}$ if and only if

$$t_{21} = t_{11}^\sigma, \ t_{22} = t_{13}^\sigma, \ t_{23} = t_{12}^\sigma, \ t_3^\sigma = t_3.$$

Note that $\nu\sigma = \sigma\nu^2$, $\nu^2\sigma = \sigma\nu$. The rest of the condition is obvious. Since $\mathrm{Gal}(k(\alpha)/k)$ is generated by the permutations (23) and (123), this proves (3).

Q.E.D.

**Theorem (4.9)** (1) *There is a bijection*

$$\mathfrak{S}_3 \setminus (k^\times / \mathrm{N}_{k_1/k}(k_1^\times))^3 / \{(t,t,t) \mid t \in k^\times\} \cong \alpha_V^{-1}(k),$$

*where the action of $\mathfrak{S}_3$ is given by permutations. Moreover $\beta = (\beta_1, \beta_2, \beta_3) \in (k^\times)^3$ corresponds to the orbit of*

$$v_1 \begin{pmatrix} \beta_1 & & \\ & -\beta_2 & \\ & & 0 \end{pmatrix} + v_2 \begin{pmatrix} 0 & & \\ & -\beta_2 & \\ & & \beta_3 \end{pmatrix}.$$

(2) *If $k'/k$ is a quadratic extension of $k$ different from $k_1$, there is a bijection*

$$\mathfrak{S}_2 \setminus k(\alpha)^\times / \mathrm{N}_{k_1(\alpha)/k(\alpha)}(k_1(\alpha)^\times) \cong \alpha_V^{-1}(k'),$$

*where $\alpha = (\alpha_1, \alpha_2, 0)$ satisfies $k(\alpha_1) = k'$, $\alpha_2$ is the conjugate of $\alpha_1$, and the action of the non-trivial element $\nu \in \mathfrak{S}_2 \cong \mathrm{Gal}(k(\alpha)/k)$ is given by the usual Galois group action. Moreover, $\beta \in k(\alpha)^\times$ corresponds to the orbit of*

$$g_\alpha \left( v_1 \begin{pmatrix} \beta & & \\ & -\beta^\nu & \\ & & 0 \end{pmatrix} + v_2 \begin{pmatrix} 0 & & \\ & -\beta^\nu & \\ & & 1 \end{pmatrix} \right).$$

(3) *The set $\alpha_V^{-1}(k_1)$ consists of a single orbit.*

(4) *If $k'/k$ is a cyclic cubic extension, there is a bijection*

$$\mathbb{Z}/3\mathbb{Z} \setminus k(\alpha_1)^\times / k^\times \mathrm{N}_{k_1(\alpha_1)/k(\alpha_1)}(k_1(\alpha_1)^\times) \cong \alpha_V^{-1}(k'),$$

*where $\alpha = (\alpha_1, \alpha_2, \alpha_3)$ satisfies $k(\alpha_1) = k'$, $\alpha_2, \alpha_3$ are the conjugates of $\alpha_1$, and the action of $\mathbb{Z}/3\mathbb{Z} \cong \mathrm{Gal}(k(\alpha)/k) = \mathrm{Gal}(k(\alpha_1)/k)$ is given by the usual Galois group action. Moreover, if $\nu \in \mathrm{Gal}(k(\alpha)/k)$ is the element satisfying $\alpha_1^\nu = \alpha_2$, $\alpha_2^\nu = \alpha_3$, $\beta \in k(\alpha_1)^\times$ corresponds to the orbit of*

$$g_\alpha \left( v_1 \begin{pmatrix} \beta & & \\ & -\beta^\nu & \\ & & 0 \end{pmatrix} + v_2 \begin{pmatrix} 0 & & \\ & -\beta^\nu & \\ & & \beta^{\nu^2} \end{pmatrix} \right).$$

(5) *If $k'/k$ is a cubic extension, whose Galois closure is an $\mathfrak{S}_3$-extension, there is a bijection*

$$k(\alpha_1)^\times / k^\times \mathrm{N}_{k_1(\alpha_1)/k(\alpha_1)}(k_1(\alpha_1)^\times) \cong \alpha_V^{-1}(k'),$$

where $\alpha, \nu$ are similar as in (4). Moreover $\beta \in k_1(\alpha_1)^\times$ corresponds to the orbit of

$$g_\alpha \left( v_1 \begin{pmatrix} \beta & -\beta^\nu \\ & 0 \end{pmatrix} + v_2 \begin{pmatrix} 0 & -\beta^\nu \\ & \beta^{\nu^2} \end{pmatrix} \right).$$

*Proof.* First note that

$$\mathrm{Aut}(\mathrm{Zero}(x))_k = \begin{cases} \mathfrak{S}_3 & k(x) = k, \\ \mathfrak{S}_2 & [k(x):k] = 2, \\ \mathbb{Z}/3\mathbb{Z} & [k(x):k] = 3, \text{ i.e., } \mathrm{Gal}(k(x)/k) = \mathbb{Z}/3\mathbb{Z}, \\ 1 & [k(x):k] = 6, \text{ i.e., } \mathrm{Gal}(k(x)/k) = \mathfrak{S}_3. \end{cases}$$

Consider (1). We have an exact sequence

$$1 \to G_w^0 \to \mathrm{GL}(1)_{k_1} \times \mathrm{GL}(1)_{k_1} \times \mathrm{GL}(1)_{k_1} \to \mathrm{GL}(1)_k \times \mathrm{GL}(1)_k \to 1,$$

where the last map is given by $(t_1, t_2, t_3) \to (\mathrm{N}_{k_1/k}(t_1 t_2^{-1}), \mathrm{N}_{k_1/k}(t_2 t_3^{-1}))$. So we have an exact sequence

$$(k_1^\times)^3 \to (k^\times)^2 \to \mathrm{H}^1(k, G_w^0) \to 1.$$

Therefore,

$$\mathrm{H}^1(k, G_w^0) \cong (k^\times/\mathrm{N}_{k_1/k}(k_1^\times))^2 \cong (k^\times/\mathrm{N}_{k_1/k}(k_1^\times))^3 / \{(t, t, t) \mid t \in k^\times\}.$$

For $\beta = (\beta_1, \beta_2) \in (k^\times)^2$, we choose a large enough finite Galois extension $k'/k$ and $\beta_1', \beta_2', \beta_3' \in k_1 \otimes k'$ so that $\mathrm{N}_{k_1/k}(\beta_1' \beta_2'^{-1}) = \beta_1$, $\mathrm{N}_{k_1/k}(\beta_2' \beta_3'^{-1}) = \beta_2$. Let

$$g_\beta = (a_3(\beta_1', \beta_2', \beta_3'), I_2).$$

Then the image of $\beta$ in $\mathrm{H}^1(k, G_w^0)$ by the boundary map is represented by $\delta g_\beta$. Clearly, $\mathfrak{S}_3$ acts by permutations of $\beta_1', \beta_2', \beta_3'$.

Consider (2). We have an exact sequence

$$1 \to G_w^0 \to \mathrm{GL}(1)_{k_1(\alpha)} \times \mathrm{GL}(1)_{k_1} \to \mathrm{GL}(1)_{k(\alpha)} \to 1,$$

where the last map is given by $(t_{11}, t_{13}) \to \mathrm{N}_{k_1(\alpha)/k(\alpha)}(t_{11} t_{13}^{-1})$. So we have an exact sequence

$$k_1(\alpha)^\times \times k_1^\times \to k(\alpha)^\times \to \mathrm{H}^1(k, G_{w_\alpha}^0) \to 1.$$

Therefore,

$$\mathrm{H}^1(k, G_{w_\alpha}^0) \cong k(\alpha)^\times / \mathrm{N}_{k_1(\alpha)/k(\alpha)}(k_1(\alpha)^\times).$$

Note that this is *not* $k(\alpha)^\times / k^\times \mathrm{N}_{k_1(\alpha)/k(\alpha)}(k_1(\alpha)^\times)$.

For $\beta \in k(\alpha)^\times$, we choose a large enough finite Galois extension $k'/k$ and $\beta' \in k_1(\alpha) \otimes k'$ so that $\mathrm{N}_{k_1/k}(\beta') = \beta$. Let

$$g_\beta = (a_3(\beta', \beta'^\nu, 1), I_2).$$

Then the image of $\beta$ in $\mathrm{H}^1(k, G^0_{w_\alpha})$ by the boundary map is represented by $\delta g_\beta$. Note that we are considering $\beta$ as a $k$–rational point of $\mathrm{GL}(1)_{k(\alpha)}$. So $\beta$ is fixed by all $\eta \in \mathrm{Gal}(k^{\mathrm{sep}}/k)$ and this is why we have $I_2$ in the definition of $g_\beta$.

The non-trivial element of $\mathrm{Aut}(\mathrm{Zero}(w_\alpha))_k$ is represented by $g_\alpha \tau_1 g_\alpha^{-1}$ and it induces the exchange of $\beta'$ and $\beta'^\nu$. Therefore, the action of $\mathrm{Aut}(\mathrm{Zero}(w_\alpha))_k \cong \mathrm{Gal}(k(\alpha)/k)$ coincides with the Galois group action.

The statement (3) follows from Lemma (1.4).

Consider (4) and (5). Let
$$A = \{t \in \mathrm{GL}(1)_{k_1(\alpha_1)} \mid \mathrm{N}_{k_1(\alpha_1)/k(\alpha_1)}(t) = 1\}.$$
For case (5), $k(\alpha) = k_1(\alpha_1)$. So for both cases, we have exact sequences

$$\begin{array}{ccccccccc}
1 & \to & A & \to & G^0_{w_\alpha} & \to & \mathrm{GL}(1)_k & \to & 1 \\
& & \| & & \downarrow & & \downarrow & & \\
1 & \to & A & \to & \mathrm{GL}(1)_{k_1(\alpha)} & \to & \mathrm{GL}(1)_{k(\alpha_1)} & \to & 1
\end{array}$$

where the last maps are given by $\mathrm{N}_{k_1(\alpha_1)/k(\alpha_1)}$.

From the above exact sequences, we have the following long exact sequences.

$$\begin{array}{ccccccccc}
\cdots & \to & k^\times & \to & \mathrm{H}^1(k, A) & \to & \mathrm{H}^1(k, G^0_{w_\alpha}) & \to & 1 \\
& & \downarrow & & \| & & & & \\
k_1(\alpha_1)^\times & \to & k(\alpha_1)^\times & \to & \mathrm{H}^1(k, A) & \to & 1 & &
\end{array}$$

Therefore, $\mathrm{H}^1(k, G^0_{w_\alpha}) \cong k(\alpha_1)^\times / k^\times \mathrm{N}_{k_1(\alpha_1)/k(\alpha_1)}(k_1(\alpha_1)^\times)$.

For $\beta \in k(\alpha_1)^\times$, we choose a large enough finite Galois extension $k'/k$ and $\beta' \in k_1(\alpha_1) \otimes k'$ so that $\mathrm{N}_{k_1(\alpha_1)/k(\alpha_1)}(\beta') = \beta$. Let

$$g_\beta = \begin{cases} (a_3(\beta', \beta'^\nu, 1), I_2) & \text{case (4)}, \\ (a_3(\beta', \beta'^\nu, \beta'^{\nu^2}), a_3(\beta'^\sigma, \beta'^{\sigma\nu}, \beta'^{\sigma\nu^2}), I_2) & \text{case (5)}. \end{cases}$$

Then the image of $\beta$ in $\mathrm{H}^1(k, G^0_{w_\alpha})$ by the boundary map is represented by $\delta g_\beta$.

In cases (4), (5), $g_\alpha \tau_2 \tau_1 g_\alpha^{-1}$ represents the element of $\mathrm{Aut}(\mathrm{Zero}(w_\alpha))_k$ which corresponds to $\nu \in \mathrm{Gal}(k(\alpha)/k)$. It maps $\beta'$ to $\beta'^\nu$. Therefore, the action of $\mathrm{Aut}(\mathrm{Zero}(w_\alpha))_k$ coincides with the Galois group action for case (4). In case (5), $g_\alpha^\sigma = g_\alpha \tau_1 \tau_2 \tau_1$ also and it maps $\beta'$ to $\beta'^\sigma$. Therefore, the action of $\mathrm{Aut}(\mathrm{Zero}(w_\alpha))_k$ coincides with the Galois group action for case (5) also.

For case (1), $g_\beta w$ is the corresponding orbit and for cases (2), (4), (5), $g_\beta w_\alpha$ is the corresponding orbit and these are the ones we stated. This proves the theorem.

Q.E.D.

**Corollary (4.10)** (1) *If $x \in V_k^{\mathrm{ss}}$, the field $k(x)$ is generated by residue fields of points in $\mathrm{Zero}(x)$.*
(2) *If $x, y \in V_k^{\mathrm{ss}}$ and $k(x) = k(y)$, $G^0_x = G^0_y$.*

*Proof.* Consider case (4) of Theorem (4.9) for example. This corollary follows from the fact that the zero set and the connected component of 1 of the stabilizer of the element

$$v_1 \begin{pmatrix} \beta & & \\ & -\beta^\nu & \\ & & 0 \end{pmatrix} + v_2 \begin{pmatrix} 0 & & \\ & -\beta^\nu & \\ & & \beta^{\nu^2} \end{pmatrix}$$

are the same as those of $w$. Other cases are similar.

Q.E.D.

## §5 Interpretation of the problems

Let $\widetilde{T} = \mathrm{Ker}(G \to \mathrm{GL}(V))$ for all the cases. If $k$ is a number field, it is possible to define the zeta function for the prehomogeneous vector space $(G/\widetilde{T}, V)$. For cases (1), (2), the convergence of the zeta function as well as the determination of the principal part is discussed in [9], [8]. For case (3), the convergence of the zeta function follows from the consideration in Part IV [11] because the weights of the representation are similar to those of the quartic case in [11].

Consider case (1). By Theorem (2.11), the orbit space parametrizes $\mathfrak{E}\mathfrak{x}_2$. Since

$$G^0_{w_\alpha}/\widetilde{T} \cong \mathrm{GL}(1)_{k_1(\alpha)}/\mathrm{GL}(1)_{k_1}$$

and $k_1$ is a fixed field, the weighting factor should be the the residue of the Dedekind zeta function for the field $k_1(\alpha)$. Therefore, we are more or less counting the class number times the regulator of fields of the form $k(\sqrt{\beta_0}, \sqrt{\beta})$ with $\beta_0$ fixed.

In cases (2) and (3) the interpretation is complicated by the fact that the map $\alpha_V$ is not injective. The expected density theorem for these cases counts rational orbits with a suitable weight and does not immediately yield a density theorem for fields since most fields are associated by $\alpha_V$ with infinitely-many orbits. The first observation to make here is that the group $G^0_x$ does not depend on $G_k x$, but only on $\alpha_V(G_k x)$ (this was verified case by case above). Thus all the orbits associated to a given field have the same weight factor and we may hope to group them together in the sum in order to obtain a density theorem for fields. Secondly in every case the group $G^0_x$ is a torus and fits into a short exact sequence whose other terms are products of groups obtained from $\mathrm{GL}(1)$ by restriction of scalars. The good behavior of Tamagawa measures in short exact sequences and under restriction of scalars (see [4]) leads us to expect that the weight factor will be essentially the class number times the regulator of the given field.

Before describing the conjectural density theorems for cases (2) and (3) it may be helpful to mention a simpler example where many of the same phenomena occur. This is the case of $G = \mathrm{GL}(2)_k$ acting on the space, $V$, of binary quadratic forms (note that, unlike [1] and [10], we are *not* including a $\mathrm{GL}(1)$ factor in $G$). Here there is a surjective map $\alpha_V : G_k \setminus V_k^{\mathrm{ss}} \to \mathfrak{E}\mathfrak{x}_2$ and if $k'$ is a quadratic extension of $k$ then

(5.1) $$\alpha_V^{-1}(k') \cong k^\times / N_{k'/k}((k')^\times).$$

It is well known (see [5] for example) that the quotient in (5.1) parametrizes cyclic algebras containing $k'$, which in this case are simply quaternion algebras. Thus, excluding the point $\alpha_V^{-1}(k)$, $G_k \setminus V_k^{\mathrm{ss}}$ may be put into one-to-one correspondence with pairs $(Q, k')$, where $Q$ is a quaternion algebra over $k$ and $k'$ is a quadratic subfield of $Q$. It is these objects which are being counted in the density theorem. In this case if $w_{k'} \in \alpha_V^{-1}(k')$ then

$$G^0_{w_{k'},k} \cong \{t \in (k')^\times \mid N_{k'/k}(t) = 1\}$$

and so the weight factor for $(Q, k')$ depends only on $k'$; it is more or less the class number times the regulator of $k'$. Different choices of $Q$ merely give orbits with different "discriminants".

Similar considerations may be applied to cases (2) and (3). In each case, by grouping together the contributions from each of the orbits in a given fiber of $\alpha_V$, we should obtain a density theorem for the class number times the regulator of a certain kind of field. For case (2) the set of fields will consist of the composita of all quadratic extensions of $k$ with a fixed cubic field. For case (3) it will consist of the composita of all fields of degree at most three with a fixed quadratic field. In case (3) the correspondence of orbits with arithmetic objects may be made one-to-one by the device of introducing cyclic algebras as in the discussion of binary quadratic forms above. For example, if $k'$ is a non-normal cubic extension, then $\alpha_V^{-1}(k')$ is

$$(k')^\times / k^\times N_{k' \cdot k_1 / k'}((k' \cdot k_1)^\times)$$

and this may be identified with the set of classes of quaternion algebras $Q$ over $k'$ which contain $k' \cdot k_1$, under the equivalence relation $Q_1 \equiv Q_2$ if there is a central simple algebra $A$ over $k$ such that $[Q_1][A \otimes k'] = [Q_2]$ in the Brauer group of $k'$. This is also possible in case (2) when $k_1/k$ is cyclic, but does not seem so easy when it is not.

One advantage of considering non-split cases in this paper is that it makes the global theory much easier. For example the group is of rank five for the split $E_6$ case in [7], and the complexity of computing the principal part of the zeta function is already formidable. However, the group is of rank three for the non-split $E_6$ case in this paper, and the global theory is well within our reach. The local theory is slightly more difficult but not much. Of course we would prefer to compute the density of the class number times the regulator of cubic fields without any modification. However, by considering the non-split $E_6$ case, we are still considering cubic fields composed with a given quadratic field and the expected density theorem will probably be reasonably satisfying.

Anthony C. Kable
Oklahoma State University
Mathematics Department
401 Math Science
Stillwater OK 74078-1058 USA
kable@math.okstate.edu

Akihiko Yukie
Oklahoma State University
Mathematics Department
401 Math Science
Stillwater OK 74078-1058 USA
yukie@math.okstate.edu